\documentclass{amsart}
\usepackage[all]{xy}
\usepackage{tikz}
 \tikzstyle{int}=[circle, draw,fill=black,outer sep=0,minimum size=3pt, inner sep=0]
  \tikzstyle{ext}=[circle, draw=black,outer sep=0,inner sep=1pt]
\usepackage{latexsym}
\usepackage{amssymb}
\usepackage{amsmath}
\usepackage{amsthm}
\usepackage{amscd}
\usepackage{fancyhdr}
\usepackage{fullpage}
\usepackage{mathrsfs}
\input cyracc.def
\xyoption{arc}
\def\id{{\mbox{1 \hskip -8pt 1}}}
\newcommand{\sgn}{{\mathit s  \mathit g\mathit  n}}
 \newcommand{\lon}{\longrightarrow}
 \newcommand{\bu}{\bullet}
 
 \newcommand{\rar}{\rightarrow}
 \newcommand{\hook}{\hookrightarrow}

\newcommand{\Id}{{\mathrm I\mathrm d}}

 \newcommand{\Z}{{\mathbb Z}}
 \newcommand{\bS}{{\mathbb S}}

 \newcommand{\R}{{\mathbb R}}
 \newcommand{\N}{{\mathbb N}}
 \newcommand{\K}{{\mathbb K}}

\newcommand{\GC}{\mathsf{GC}}
\newcommand{\fGC}{\mathsf{fGC}}

\newcommand{\wh}{\widehat}

 \newcommand{\ot}{\otimes}

\newcommand{\sG}{{\mathsf G}}

\newcommand{\Def}{{\mathsf D\mathsf e\mathsf f }}

 \newcommand{\Beq}{\begin{equation}}
 \newcommand{\Eeq}{\end{equation}}
 \newcommand{\Beqr}{\begin{eqnarray}}
 \newcommand{\Eeqr}{\end{eqnarray}}
 \newcommand{\Beqrn}{\begin{eqnarray*}}
 \newcommand{\Eeqrn}{\end{eqnarray*}}
 \newcommand{\Ba}{\begin{array}}
 \newcommand{\Ea}{\end{array}}
 \newcommand{\Bi}{\begin{itemize}}
 \newcommand{\Ei}{\end{itemize}}
 \newcommand{\Bc}{\begin{center}}
 \newcommand{\Ec}{\end{center}}

 \newcommand{\f}{{\mathcal O}}
 \newcommand{\cA}{{\mathcal A}}

 \newcommand{\caD}{{\mathcal D}}
 \newcommand{\cE}{{\mathcal E}}
 \newcommand{\cF}{{\mathcal F}}

 \newcommand{\caL}{{\mathcal L}}
 \newcommand{\cM}{{\mathcal M}}
 
 \newcommand{\cP}{{\mathcal P}}
 \newcommand{\cQ}{{\mathcal Q}}
 \newcommand{\cR}{{\mathcal R}}
 
 \newcommand{\cT}{{\mathcal T}}

 \newcommand{\al}{\alpha}
 \newcommand{\be}{\beta}
 \newcommand{\ga}{\gamma}
 
 \newcommand{\Ga}{\Gamma}

 \newcommand{\la}{\lambda}

 \newcommand{\Hom}{{\mathrm H\mathrm o\mathrm m}}
 
 \newcommand{\sip}{\smallskip}
 \newcommand{\bip}{\bigskip}
 \newcommand{\mip}{\vspace{2.5mm}}

  \newcommand{\dGC}{\mathsf{dGC}}
   \newcommand{\dcGC}{\mathsf{dcGC}}

   \newcommand{\Ass}{\mathcal{A} \mathit{ss}}
 \newcommand{\swhHoLBcd}{\mathcal{H}\mathit{olieb}_{c,d}^{_{\atop ^{\bigstar\circlearrowright}}}}
  \newcommand{\swHoLBcd}{\swhHoLBcd}
  \newcommand{\swhHoLB}{\mathcal{H}\mathit{olieb}^{_{\atop ^{\bigstar\circlearrowright}}}}
\newcommand{\swHoLB}{\swhHoLB}
 \newcommand{\sHoLBcd}{\mathcal{H}\mathit{olieb}_{c,d}^{_{\atop ^{\bigstar}}}}
 \newcommand{\sHoLB}{\mathcal{H}\mathit{olieb}^{_{\atop ^{\bigstar}}}}

 \newcommand{\suHoLBcd}{\mathcal{H}\mathit{olieb}_{c,d}^{_{\atop ^{\bigstar\uparrow}}}}

  \newcommand{\LB}{\mathcal{L}\mathit{ieb}}

\newcommand{\LBcd}{\mathcal{L}\mathit{ieb}_{c,d}}

\newcommand{\HoLBcd}{\mathcal{H}\mathit{olieb}_{c,d}}

\newcommand{\wHoLBcd}{\widehat{\mathcal{H}\mathit{olieb}}_{c,d}}

\newcommand{\HoLB}{\mathcal{H}\mathit{olieb}}

 \newcommand{\grt}{\mathfrak{grt}}
 \newcommand{\Der}{\mathrm{Der}}

\theoremstyle{plain}
\swapnumbers

\newtheorem{prop-def}[theorem]{Proposition-definition}

\newtheorem{main-theorem}{Main~Theorem}[section]
\newtheorem{section-theorem}{Theorem}[section]
\newtheorem{section-corollary}{Corollary}[section]

\theoremstyle{definition}

\begin{document}

\sloppy

 \newenvironment{proo}{\begin{trivlist} \item{\sc {Proof.}}}
  {\hfill $\square$ \end{trivlist}}

\long\def\symbolfootnote[#1]#2{\begingroup%
\def\thefootnote{\fnsymbol{footnote}}\footnote[#1]{#2}\endgroup}

 \title{From deformation theory of wheeled props to\\ classification of Kontsevich formality maps}

\author{Assar~Andersson}
\address{Assar~Andersson:  Mathematics Research Unit, Luxembourg University, Maison du Nombre, 6 Avenue de la Fonte,
 L-4364 Esch-sur-Alzette,   Grand Duchy of Luxembourg }
\email{assar.andersson@uni.lu}
\author{Sergei~Merkulov}
\address{Sergei~Merkulov:  Mathematics Research Unit, Luxembourg University, Maison du Nombre, 6 Avenue de la Fonte,
 L-4364 Esch-sur-Alzette,   Grand Duchy of Luxembourg }
\email{sergei.merkulov@uni.lu}


 \begin{abstract}
We study the homotopy theory of the wheeled prop controlling Poisson structures on formal graded finite-dimensional manifolds and prove, in particular, that
Grothendieck-Teichm\"uller group acts on that wheeled prop faithfully and homotopy non-trivially. Next we apply this homotopy theory to the study of the deformation complex of an arbitrary
M.\ Kontsevich formality map and compute the full cohomology group of that deformation complex in terms of the cohomology of a certain graph complex introduced earlier by M.\ Kontsevich in \cite{Ko} and studied by T.\ Willwacher in \cite{Wi1}.

\end{abstract}
 \maketitle

{\large
\section{\bf Introduction}
}
\label{sec:introduction}

\subsection{Wheeled props, formal Poisson structures and Grothendieck-Teichm\"uller group}
Let $V$ be an arbitrary finite-dimensional $\Z$-graded vector space over a field $\K$ of characteristic zero (say, $V=\R^d$, $\K=\R$) and $V^*:=\Hom(V,\K)$ its dual.
Then the completed symmetric algebra $\f_\cM:=\widehat{\odot^\bu}V$ can be understood as the $\K$-algebra of formal smooth functions on the dual vector space $V^*$ understood as a formal manifold $\cM$, and the Lie algebra of derivations of $\f_\cM$,
$$
\cT\cM:= \Der(\f_V)\simeq \Hom(V,\widehat{\odot^\bu} V)\simeq \prod_{m\geq 0}
\Hom(V, \odot^m V),
$$
as the Lie algebra of formal smooth vector fields on $\cM$. A {\em formal graded Poisson structure}\,
on $\cM$ is a degree 2 element $\pi$ in the Schouten Lie algebra
$$
\cT_{poly}{\cM}:=\wedge^\bu T_\cM \simeq \prod_{m\geq 0,n\geq 0}
\Hom(\wedge^n V, \odot^m V)[-n]=\prod_{m\geq 0,n\geq 0}
\Hom(\odot^n (V[1]), \odot^m V)=\prod_{k\geq 0}\odot^k (V^*[-1]\oplus V)
$$
of polyvector fields satisfying the Maurer-Cartan equation,
$$
[\pi,\pi]=0,
$$
where the Schouten Lie bracket $[\ ,\ ]$ (of degree $-1$) originates essentially from the canonical pairing map $V^*[-1]\times V \rar \K[-1]$. Thus a formal Poisson structure is a formal power series,
\Beq\label{1: pi formal power series}
\pi=\sum_{n,m=0}^\infty\pi_n^m,\ \ \ \  \ \pi_n^m\in \Hom(\wedge^n V, \odot^m V)[2-n]
\Eeq
and hence can be understood as a representation in the vector space $V$,
$$
\rho_\pi: \HoLB_{0,1}^{\star} \lon \cE nd_V,
$$
of a certain prop of {\em formal Poisson structures}\, which is by definition the free
prop\footnote{See, e.g., \cite{Ma,Me3} and the first sections of \cite{V}  for an elementary introduction into the theory of props and wheeled props.} generated a collection of 1-dimensional $\bS_m^{op}\times \bS_n$ bimodules $\id_m\ot \sgn_n[n-2]$. A useful observation \cite{Me0} is that this prop comes equipped with a natural differential $\delta^\star$ such that the Maurer-Cartan equation $[\pi,\pi]=0$ gets encoded into the compatibility of the representation $\rho_\pi$ with the differentials. The strange notation $\HoLB_{1,0}^{\star}$ comes from the fact that this particular prop comes from the family of dg props $\HoLB_{c,d}^{\star}$ which control Maurer-Cartan elements $\pi$ in the graded commutative algebra
\Beq\label{1: odot V[-c] + V[-d]}
\prod_{k\geq 0}\odot^k (V^*[-d]\oplus V[-c])
\Eeq
equipped with the obvious Poisson type Lie bracket (of homological degree $-c-d$). The case $c=0$, $d=1$ corresponds to formal Poisson structures while the case $c=1$, $d=1$ corresponds to extended {\em homotopy Lie bialgebras}\, structures.

\sip

The superscript $\star$ in the notation indicates that we consider in this paper an {\em extended
version}\,  of the family of props $\HoLBcd$ studied earlier in \cite{MW1}. The latter family controls the {\em truncated}\,  version of the above formal power series (\ref{1: pi formal power series}) which allows only monomials with $m,n\geq 1$, $m+n\geq 3$; such a truncation makes perfect sense in the context of the theory of minimal resolutions of $(c,d)$ Lie bialgebras. However we are interested in this paper in the ``full story" with no restrictions
on the integer parameters $m$ and $n$, and that ``full story" turns out to be sometimes quite different from the truncated one.

\sip

Note that under appropriate completion of the above graded commutative algebra (\ref{1: odot V[-c] + V[-d]}) all the above
structures (the convergent Lie bracket, Maurer-Cartan elements $\pi$ and associated representations $\rho_\pi$ of the props $\HoLBcd^\star$) make sense also for {\em infinite-dimensional}\, vector spaces $V$. An important point of this paper is that the deformation theories of these structures behave
quite differently in finite and infinite dimensions. Indeed, in infinite dimensions they can be understood as representations of {\em ordinary}\, props
$$
\rho_\pi: \sHoLBcd \lon \cE nd_V
$$
while in {\em finite}\, dimensions as representations of their {\em wheeled closures} (cf.\ \cite{Me1,MMS}),
$$
\rho_\pi: \swHoLBcd \lon \cE nd_V
$$
which have quite different deformation theories or, equivalently, quite different dg Lie algebras,
$\Der(\sHoLBcd)$ and $\Der(\swHoLBcd)$, of derivations (see \S 3 for their precise definitions).
 The dg prop $\sHoLBcd$ is a proper subprop of $\swHoLBcd$, the latter containing many more universal operations (involving, roughly speaking, the trace operation $V\ot V^*\rar \K$ which has no sense in general when $\dim V=\infty$).
\sip

The first main purpose of this paper is the study of the deformation theory of both props
$\sHoLBcd$ and $\swHoLBcd$ (in fact of their completed versions) and the computation of the cohomologies of the associated complexes of derivations  in terms of the M.\ Kontsevich graph complexes
$\GC_d^{\geq 2}$ introduced\footnote{The symbol $\GC_d$ stands often in the literature for the graph complex generated by connected oriented graphs with all vertices trivalent; we denote by $\GC_d^{\geq 2}$ its extension which allows connected graphs with at least bivalent vertices (see \S 3.2 for more details and references). The latter complex has a quasi-isomorphic versions, $\dGC_d^{\geq 2}$, spanned by graphs with fixed directions on edges; the subcomplex of $\dGC_d^{\geq 2}$ spanned by {\em oriented}\, graphs, that is, directed graphs with no closed paths of directed edges, is denoted by $\GC_d^{or}$. These complexes have been studied in \cite{Wi1,Wi2,Z}.}  in \cite{Ko} and studied in \cite{Wi1}, and of its oriented version $\GC_d^{or}$
which was studied in \cite{Wi2,Z}. These complexes are spanned by {\em connected graphs}. It is often useful  \cite{MW1,MW3} to add to these classical graph complexes an additional element $\emptyset$ concentrated in degree zero, ``a graph with no vertices and edges", and define the {\em full graph complexes}\, of not necessarily connected graphs as the completed graded symmetric tensor algebras
\Beq\label{1: fGC in terms of GC}
 \fGC^{\geq 2}_{d}:=\widehat{\odot^\bu}\left((\GC_{d}^{\geq 2}\oplus \K )[-d]\right)[d],
\Eeq
\Beq\label{1: fGC^or in terms of GC^0r}
\fGC^{or}_{d}:=\widehat{\odot^\bu}\left((\GC_{d}^{or}\oplus \K )[-1-c-d]\right)[d],
\Eeq
the summands $\K$ being generated by $\emptyset$. The formal class $\emptyset$ takes care
for (homotopy non-trivial) rescaling operations of the (wheeled) props under considerations, and essentially leads us in applications to the full Grothendieck-Teichm\"uller group $GRT=GRT_1  \rtimes\K^*$ (see \cite{D2}) rather than to its reduced version $GRT_1$. The Lie bracket of $\emptyset$  with elements $\Ga$  of $\GC_{d}^{\geq 2}$ or $\GC_{d}^{or}$ is defined as the multiplication of $\Ga$  by twice the number of its loops.

\subsubsection{\bf Proposition} {\em There are morphisms of dg Lie algebras,
$$
F^\circlearrowright: \fGC^{\geq 2}_{c+d+1} \lon \Der(\swHoLBcd), \ \ \ \ F: \fGC_{c+d+1}^{or}
\lon  \Der(\sHoLBcd)
$$
which are quasi-isomorphisms.
}


\sip

It was proven in \cite{Wi1,Wi2} that $H^\bu(\GC_{c+d+1}^{\geq 2}) =H^\bu(\GC_{c+d+2}^{or})$ and that
$$
H^0(\fGC_2^{\geq 2})=H^0(\fGC_3^{or})=\grt\ \stackrel{\text{as a vector space}}{\simeq}\ \grt_1 \oplus \K,
$$
where $\grt$ (resp., $\grt_1$) is the Lie algebra of the Grothendieck-Teichm\"uller group $GRT$ (resp., $GRT_1$). It is easy to see that
$H^0(\GC_2^{or})=0$ and $H^0(\GC_3^{\geq 2})=0$.

\subsubsection{\bf Corollary} {\em There is an isomorphism of Lie algebras
$$
H^0(\Der(\swHoLB_{0,1}))=\grt
$$
that is, the Grothendieck-Teichm\"uller group $GRT$ acts up to homotopy faithfully (and essentially transitively) on the  vertex completion of the wheeled properad $\swHoLB_{1,0}$ governing {\em finite}-dimensional formal Poisson structures.

\sip

By contrast
$$
H^0(\Der(\sHoLB_{0,1}))=0
$$
that is,
the completion of the properad $\sHoLB_{1,0}$ governing {\em infinite}-dimensional formal Poisson structures admits no homotopy non-trivial automorphisms at all.
}

\sip

Note that the above Proposition applied to another interesting case $c=d=1$ gives us quite the opposite picture,
$$
H^0(\Der(\swHoLB_{1,1}))=0, \ \ \ \ \ H^0(\Der(\sHoLB_{1,1}))=\grt
$$
cf.\ \cite{MW1}. These results are by no means surprising --- the graph complex
$\fGC_2^{\geq 2}$ can be understood as a kind of universal incarnation of the Chevalley-Eilenberg deformation complex of the Lie algebra $\cT_{poly}(\cM)$ for any finite-dimensional formal manifold \cite{Ko}, and the fact that $H^0(\fGC_2^{\geq 2})=\grt$
already implies \cite{Wi1} that the Grothendieck-Teichm\"uller group $GRT$ acts (up to homotopy) as universal $\caL ie_\infty$ automorphisms of $\cT_{poly}(\cM)$; this action is given in terms of certain iterations of the canonical $GL(V)$-invariant $BV$ operator on  $\cT_{poly}(\cM)$, so what the above Corollary says essentially is that even if one drops this restriction on the possible structure of linear operators acting on $\cT_{poly}(\cM)$, the action of $GRT$ remains homotopy non-trivial.
 The above Proposition can be inferred from the theory of stable cohomology of the Lie algebra of polyvector fields developed in \cite{Wi3} (but not immediately).
In any case, our proof of Proposition 1.1.1 is very short, so we decided to show a new direct argument behind that claim in \S 4.1 below.

\sip

 The main advantage of our study of the homotopy theory of the vertex completion $\wh{\HoLB}_{0,1}^{\atop \bigstar\circlearrowright}$ of the wheeled prop $\swHoLB_{0,1}$ is that it gives us --- almost immediately! --- a full insight into the homotopy theory of M.\ Kontsevich's formality maps which is the second main topic of this paper.

 \subsection{Homotopy classification of M.\ Kontsevich formality maps} M.\ Kontsevich formality map \cite{Ko2} associates to any {\em finite}-dimensional formal Poisson structure $\pi$ on a formal graded manifold $\cM=V^*$ a curved $\cA ss_\infty$-structure on the $\R$-algebra $\f_\cM=\widehat{\odot^\bu} V$ of formal smooth functions on $\cM$ which is given in terms of polydifferential operators constructed from $\pi$. In our approach $\pi$ is a representation in $V$ of the wheeled prop $\swHoLB_{1,0}$, and the construction of  polydifferentials operators
 from $\pi$ can be conveniently encoded into the polydifferential functor \cite{MW3}
 $$
 \f: \mathsf{Category\ of\ dg\ props} \lon \mathsf{Category\ of\ dg\ operads}
 $$
 applied to the prop $\swHoLB_{1,0}$: for any dg prop $\cP$ the associated dg operad $\f(\cP)$ has the property that for any representation $\rho$ of $\cP$  in a vector space $V$ the operad $\f(\cP)$ has a canonically associated representation $\f(\rho)$ in the completed  graded commutative algebra $\widehat{\odot^\bu} V$ given in terms of polydifferenial operators. Curved $\cA ss_\infty$
algebra structures are controlled by the well-known (non-cofibrant) dg operad $c\cA ss_\infty$ so that the Maxim Kontsevich universal formality map from \cite{Ko} (or any other universal formality map) can be understood as a morphism of dg operads
$$
\cF:{c \cA ss}_\infty \lon  \f(\wh{\HoLB}_{0,1}^{\atop \bigstar\circlearrowright}).
$$
satisfying certain non-triviality conditions (se \S 5 for details). We show in this paper a very short and elementary (based essentially on the contractility of the permutahedra polytopes) proof of the following classification theorem.

\subsubsection{\bf Theorem} {\em Let\, $\Def\left(c\cA ss_\infty \stackrel{\cF}{\rar}  \f(\swHoLB_{1,0})\right)$ be the deformation complex of any given  formality map $\cF$ (in particular, of the M.\ Kontsevich map from \cite{Ko2}). Then there is a canonical morphism of complexes
 $$
 \fGC_2^{\geq 2} \lon \Def\left(c\cA ss_\infty \stackrel{\cF}{\rar} \f(\swHoLB_{0,1})\right)[1]
 $$
 which is a quasi-isomorphism.}

  \mip

This result implies the equality of cohomology groups for any $i\in \Z$,
$$
H^{i+1}\left(\Def\left(c\cA ss_\infty \stackrel{\cF}{\rar} \f(\wh{\HoLB}_{0,1}^{\atop \bigstar\circlearrowright})\right)\right)=H^i(\fGC_2^{\geq 2})
$$
which in the special case $i=0$ reads as
$$
H^1\left(\Def\left(c\cA ss_\infty \stackrel{\cF}{\rar} \f(\wh{\HoLB}_{0,1}^{\atop \bigstar\circlearrowright})\right)\right)=H^0(\fGC_2^{\geq 2}) =\grt
$$
and hence gives us a new (very short) proof of the following remarkable Theorem by V.\ Dolgushev.

\subsubsection{\bf Theorem \cite{Do}} {\em The Grothendieck-Teichm\"uller  group $GRT$  acts freely and
transitively on the set of homotopy classes of universal formality morphisms.}

\sip

This Theorem implies the identification of the set of homotopy classes of formality maps
with the set of V.\ Drinfeld associators \cite{D2}.

\subsection{Some notation}
 The set $\{1,2, \ldots, n\}$ is abbreviated to $[n]$;  the group of bijections $[n]\rar [n]$ is
denoted by $\bS_n$;
the trivial (resp., sign) one-dimensional representation of
 $\bS_n$ is denoted by $\id_n$ (resp., $\sgn_n$).
The cardinality of a finite set $S$ is denoted by $\# S$.
We work in this paper in the category of $\Z$-graded vector spaces over a field $\K$
of characteristic zero.
If $V=\oplus_{i\in \Z} V^i$ is a graded vector space, then
$V[k]$ stands for the graded vector space with $V[k]^i:=V^{i+k}$; for $v\in V^i$ we set $|v|:=i$.
If $V$ is a complex with a differential $d$, then $V[k]$ is also a complex with the differential given by $(-1)^k d$.

\sip

For the basic notions and  facts of the theory of props and properads we refer to the papers \cite{Ma,MV,V} (and references cited there) and of their wheeled versions to \cite{MMS, Me1}. A short introduction into these theories can be found in  \cite{Me3}.
We assume that every (wheeled) prop $\cP$ we work with in this paper has the unit denote by $\uparrow\in \cP(1,1)$.

\sip

\bip

{
\Large
\section{\bf Wheeled properads of homotopy Lie bialgebras and their extensions}
}

\sip

\subsection{Reminder on props of  Lie $(c,d)$-bialgebras and their minimal resolutions} Consider for any pair of integeres $c,d\in \Z$ a quadratic prop \cite{MW1}
$$
\LB_{c,d}:=\cF ree\langle e\rangle/\langle\cR\rangle,
$$
defined as the quotient
of the free prop generated by an  $\bS$-bimodule $e=\{e(m,n)\}_{m,n\geq 0}$ with
 all $e(m,n)=0$ except\footnote{When representing elements of various props
  below as graphs we always assume by default that all edges and legs are {\em directed}\, with the flow running from the bottom of the graph to the top.}
$$
e(2,1):=\id_1\ot \sgn_2^{c}[c-1]=\mbox{span}\left\langle
\Ba{c}\begin{xy}
 <0mm,-0.55mm>*{};<0mm,-2.5mm>*{}**@{-},
 <0.5mm,0.5mm>*{};<2.2mm,2.2mm>*{}**@{-},
 <-0.48mm,0.48mm>*{};<-2.2mm,2.2mm>*{}**@{-},
 <0mm,0mm>*{\circ};<0mm,0mm>*{}**@{},
 <0.5mm,0.5mm>*{};<2.7mm,2.8mm>*{^{_2}}**@{},
 <-0.48mm,0.48mm>*{};<-2.7mm,2.8mm>*{^{_1}}**@{},
 \end{xy}\Ea
=(-1)^{c}
\Ba{c}\begin{xy}
 <0mm,-0.55mm>*{};<0mm,-2.5mm>*{}**@{-},
 <0.5mm,0.5mm>*{};<2.2mm,2.2mm>*{}**@{-},
 <-0.48mm,0.48mm>*{};<-2.2mm,2.2mm>*{}**@{-},
 <0mm,0mm>*{\circ};<0mm,0mm>*{}**@{},
 <0.5mm,0.5mm>*{};<2.7mm,2.8mm>*{^{_1}}**@{},
 <-0.48mm,0.48mm>*{};<-2.7mm,2.8mm>*{^{_2}}**@{},
 \end{xy}\Ea
   \right\rangle
$$
$$
e(1,2):= \sgn_2^{d}\ot \id_1[d-1]=\mbox{span}\left\langle
\Ba{c}\begin{xy}
 <0mm,0.66mm>*{};<0mm,3mm>*{}**@{-},
 <0.39mm,-0.39mm>*{};<2.2mm,-2.2mm>*{}**@{-},
 <-0.35mm,-0.35mm>*{};<-2.2mm,-2.2mm>*{}**@{-},
 <0mm,0mm>*{\circ};<0mm,0mm>*{}**@{},
   <0.39mm,-0.39mm>*{};<2.9mm,-4mm>*{^{_2}}**@{},
   <-0.35mm,-0.35mm>*{};<-2.8mm,-4mm>*{^{_1}}**@{},
\end{xy}\Ea
=(-1)^{d}
\Ba{c}\begin{xy}
 <0mm,0.66mm>*{};<0mm,3mm>*{}**@{-},
 <0.39mm,-0.39mm>*{};<2.2mm,-2.2mm>*{}**@{-},
 <-0.35mm,-0.35mm>*{};<-2.2mm,-2.2mm>*{}**@{-},
 <0mm,0mm>*{\circ};<0mm,0mm>*{}**@{},
   <0.39mm,-0.39mm>*{};<2.9mm,-4mm>*{^{_1}}**@{},
   <-0.35mm,-0.35mm>*{};<-2.8mm,-4mm>*{^{_2}}**@{},
\end{xy}\Ea
\right\rangle
$$
by the ideal generated by the following elements
\Beq\label{R for LieB}
\cR:\left\{
\Ba{c}
\Ba{c}\resizebox{7mm}{!}{
\begin{xy}
 <0mm,0mm>*{\circ};<0mm,0mm>*{}**@{},
 <0mm,-0.49mm>*{};<0mm,-3.0mm>*{}**@{-},
 <0.49mm,0.49mm>*{};<1.9mm,1.9mm>*{}**@{-},
 <-0.5mm,0.5mm>*{};<-1.9mm,1.9mm>*{}**@{-},
 <-2.3mm,2.3mm>*{\circ};<-2.3mm,2.3mm>*{}**@{},
 <-1.8mm,2.8mm>*{};<0mm,4.9mm>*{}**@{-},
 <-2.8mm,2.9mm>*{};<-4.6mm,4.9mm>*{}**@{-},
   <0.49mm,0.49mm>*{};<2.7mm,2.3mm>*{^3}**@{},
   <-1.8mm,2.8mm>*{};<0.4mm,5.3mm>*{^2}**@{},
   <-2.8mm,2.9mm>*{};<-5.1mm,5.3mm>*{^1}**@{},
 \end{xy}}\Ea
 +
\Ba{c}\resizebox{7mm}{!}{\begin{xy}
 <0mm,0mm>*{\circ};<0mm,0mm>*{}**@{},
 <0mm,-0.49mm>*{};<0mm,-3.0mm>*{}**@{-},
 <0.49mm,0.49mm>*{};<1.9mm,1.9mm>*{}**@{-},
 <-0.5mm,0.5mm>*{};<-1.9mm,1.9mm>*{}**@{-},
 <-2.3mm,2.3mm>*{\circ};<-2.3mm,2.3mm>*{}**@{},
 <-1.8mm,2.8mm>*{};<0mm,4.9mm>*{}**@{-},
 <-2.8mm,2.9mm>*{};<-4.6mm,4.9mm>*{}**@{-},
   <0.49mm,0.49mm>*{};<2.7mm,2.3mm>*{^2}**@{},
   <-1.8mm,2.8mm>*{};<0.4mm,5.3mm>*{^1}**@{},
   <-2.8mm,2.9mm>*{};<-5.1mm,5.3mm>*{^3}**@{},
 \end{xy}}\Ea
 +
\Ba{c}\resizebox{7mm}{!}{\begin{xy}
 <0mm,0mm>*{\circ};<0mm,0mm>*{}**@{},
 <0mm,-0.49mm>*{};<0mm,-3.0mm>*{}**@{-},
 <0.49mm,0.49mm>*{};<1.9mm,1.9mm>*{}**@{-},
 <-0.5mm,0.5mm>*{};<-1.9mm,1.9mm>*{}**@{-},
 <-2.3mm,2.3mm>*{\circ};<-2.3mm,2.3mm>*{}**@{},
 <-1.8mm,2.8mm>*{};<0mm,4.9mm>*{}**@{-},
 <-2.8mm,2.9mm>*{};<-4.6mm,4.9mm>*{}**@{-},
   <0.49mm,0.49mm>*{};<2.7mm,2.3mm>*{^1}**@{},
   <-1.8mm,2.8mm>*{};<0.4mm,5.3mm>*{^3}**@{},
   <-2.8mm,2.9mm>*{};<-5.1mm,5.3mm>*{^2}**@{},
 \end{xy}}\Ea
 \ \ , \ \
\Ba{c}\resizebox{8.4mm}{!}{ \begin{xy}
 <0mm,0mm>*{\circ};<0mm,0mm>*{}**@{},
 <0mm,0.69mm>*{};<0mm,3.0mm>*{}**@{-},
 <0.39mm,-0.39mm>*{};<2.4mm,-2.4mm>*{}**@{-},
 <-0.35mm,-0.35mm>*{};<-1.9mm,-1.9mm>*{}**@{-},
 <-2.4mm,-2.4mm>*{\circ};<-2.4mm,-2.4mm>*{}**@{},
 <-2.0mm,-2.8mm>*{};<0mm,-4.9mm>*{}**@{-},
 <-2.8mm,-2.9mm>*{};<-4.7mm,-4.9mm>*{}**@{-},
    <0.39mm,-0.39mm>*{};<3.3mm,-4.0mm>*{^3}**@{},
    <-2.0mm,-2.8mm>*{};<0.5mm,-6.7mm>*{^2}**@{},
    <-2.8mm,-2.9mm>*{};<-5.2mm,-6.7mm>*{^1}**@{},
 \end{xy}}\Ea
 +
\Ba{c}\resizebox{8.4mm}{!}{ \begin{xy}
 <0mm,0mm>*{\circ};<0mm,0mm>*{}**@{},
 <0mm,0.69mm>*{};<0mm,3.0mm>*{}**@{-},
 <0.39mm,-0.39mm>*{};<2.4mm,-2.4mm>*{}**@{-},
 <-0.35mm,-0.35mm>*{};<-1.9mm,-1.9mm>*{}**@{-},
 <-2.4mm,-2.4mm>*{\circ};<-2.4mm,-2.4mm>*{}**@{},
 <-2.0mm,-2.8mm>*{};<0mm,-4.9mm>*{}**@{-},
 <-2.8mm,-2.9mm>*{};<-4.7mm,-4.9mm>*{}**@{-},
    <0.39mm,-0.39mm>*{};<3.3mm,-4.0mm>*{^2}**@{},
    <-2.0mm,-2.8mm>*{};<0.5mm,-6.7mm>*{^1}**@{},
    <-2.8mm,-2.9mm>*{};<-5.2mm,-6.7mm>*{^3}**@{},
 \end{xy}}\Ea
 +
\Ba{c}\resizebox{8.4mm}{!}{ \begin{xy}
 <0mm,0mm>*{\circ};<0mm,0mm>*{}**@{},
 <0mm,0.69mm>*{};<0mm,3.0mm>*{}**@{-},
 <0.39mm,-0.39mm>*{};<2.4mm,-2.4mm>*{}**@{-},
 <-0.35mm,-0.35mm>*{};<-1.9mm,-1.9mm>*{}**@{-},
 <-2.4mm,-2.4mm>*{\circ};<-2.4mm,-2.4mm>*{}**@{},
 <-2.0mm,-2.8mm>*{};<0mm,-4.9mm>*{}**@{-},
 <-2.8mm,-2.9mm>*{};<-4.7mm,-4.9mm>*{}**@{-},
    <0.39mm,-0.39mm>*{};<3.3mm,-4.0mm>*{^1}**@{},
    <-2.0mm,-2.8mm>*{};<0.5mm,-6.7mm>*{^3}**@{},
    <-2.8mm,-2.9mm>*{};<-5.2mm,-6.7mm>*{^2}**@{},
 \end{xy}}\Ea
 \\
 \Ba{c}\resizebox{5mm}{!}{\begin{xy}
 <0mm,2.47mm>*{};<0mm,0.12mm>*{}**@{-},
 <0.5mm,3.5mm>*{};<2.2mm,5.2mm>*{}**@{-},
 <-0.48mm,3.48mm>*{};<-2.2mm,5.2mm>*{}**@{-},
 <0mm,3mm>*{\circ};<0mm,3mm>*{}**@{},
  <0mm,-0.8mm>*{\circ};<0mm,-0.8mm>*{}**@{},
<-0.39mm,-1.2mm>*{};<-2.2mm,-3.5mm>*{}**@{-},
 <0.39mm,-1.2mm>*{};<2.2mm,-3.5mm>*{}**@{-},
     <0.5mm,3.5mm>*{};<2.8mm,5.7mm>*{^2}**@{},
     <-0.48mm,3.48mm>*{};<-2.8mm,5.7mm>*{^1}**@{},
   <0mm,-0.8mm>*{};<-2.7mm,-5.2mm>*{^1}**@{},
   <0mm,-0.8mm>*{};<2.7mm,-5.2mm>*{^2}**@{},
\end{xy}}\Ea
  -
\Ba{c}\resizebox{7mm}{!}{\begin{xy}
 <0mm,-1.3mm>*{};<0mm,-3.5mm>*{}**@{-},
 <0.38mm,-0.2mm>*{};<2.0mm,2.0mm>*{}**@{-},
 <-0.38mm,-0.2mm>*{};<-2.2mm,2.2mm>*{}**@{-},
<0mm,-0.8mm>*{\circ};<0mm,0.8mm>*{}**@{},
 <2.4mm,2.4mm>*{\circ};<2.4mm,2.4mm>*{}**@{},
 <2.77mm,2.0mm>*{};<4.4mm,-0.8mm>*{}**@{-},
 <2.4mm,3mm>*{};<2.4mm,5.2mm>*{}**@{-},
     <0mm,-1.3mm>*{};<0mm,-5.3mm>*{^1}**@{},
     <2.5mm,2.3mm>*{};<5.1mm,-2.6mm>*{^2}**@{},
    <2.4mm,2.5mm>*{};<2.4mm,5.7mm>*{^2}**@{},
    <-0.38mm,-0.2mm>*{};<-2.8mm,2.5mm>*{^1}**@{},
    \end{xy}}\Ea
  - (-1)^{d}
\Ba{c}\resizebox{7mm}{!}{\begin{xy}
 <0mm,-1.3mm>*{};<0mm,-3.5mm>*{}**@{-},
 <0.38mm,-0.2mm>*{};<2.0mm,2.0mm>*{}**@{-},
 <-0.38mm,-0.2mm>*{};<-2.2mm,2.2mm>*{}**@{-},
<0mm,-0.8mm>*{\circ};<0mm,0.8mm>*{}**@{},
 <2.4mm,2.4mm>*{\circ};<2.4mm,2.4mm>*{}**@{},
 <2.77mm,2.0mm>*{};<4.4mm,-0.8mm>*{}**@{-},
 <2.4mm,3mm>*{};<2.4mm,5.2mm>*{}**@{-},
     <0mm,-1.3mm>*{};<0mm,-5.3mm>*{^2}**@{},
     <2.5mm,2.3mm>*{};<5.1mm,-2.6mm>*{^1}**@{},
    <2.4mm,2.5mm>*{};<2.4mm,5.7mm>*{^2}**@{},
    <-0.38mm,-0.2mm>*{};<-2.8mm,2.5mm>*{^1}**@{},
    \end{xy}}\Ea
  - (-1)^{d+c}
\Ba{c}\resizebox{7mm}{!}{\begin{xy}
 <0mm,-1.3mm>*{};<0mm,-3.5mm>*{}**@{-},
 <0.38mm,-0.2mm>*{};<2.0mm,2.0mm>*{}**@{-},
 <-0.38mm,-0.2mm>*{};<-2.2mm,2.2mm>*{}**@{-},
<0mm,-0.8mm>*{\circ};<0mm,0.8mm>*{}**@{},
 <2.4mm,2.4mm>*{\circ};<2.4mm,2.4mm>*{}**@{},
 <2.77mm,2.0mm>*{};<4.4mm,-0.8mm>*{}**@{-},
 <2.4mm,3mm>*{};<2.4mm,5.2mm>*{}**@{-},
     <0mm,-1.3mm>*{};<0mm,-5.3mm>*{^2}**@{},
     <2.5mm,2.3mm>*{};<5.1mm,-2.6mm>*{^1}**@{},
    <2.4mm,2.5mm>*{};<2.4mm,5.7mm>*{^1}**@{},
    <-0.38mm,-0.2mm>*{};<-2.8mm,2.5mm>*{^2}**@{},
    \end{xy}}\Ea
 - (-1)^{c}
\Ba{c}\resizebox{7mm}{!}{\begin{xy}
 <0mm,-1.3mm>*{};<0mm,-3.5mm>*{}**@{-},
 <0.38mm,-0.2mm>*{};<2.0mm,2.0mm>*{}**@{-},
 <-0.38mm,-0.2mm>*{};<-2.2mm,2.2mm>*{}**@{-},
<0mm,-0.8mm>*{\circ};<0mm,0.8mm>*{}**@{},
 <2.4mm,2.4mm>*{\circ};<2.4mm,2.4mm>*{}**@{},
 <2.77mm,2.0mm>*{};<4.4mm,-0.8mm>*{}**@{-},
 <2.4mm,3mm>*{};<2.4mm,5.2mm>*{}**@{-},
     <0mm,-1.3mm>*{};<0mm,-5.3mm>*{^1}**@{},
     <2.5mm,2.3mm>*{};<5.1mm,-2.6mm>*{^2}**@{},
    <2.4mm,2.5mm>*{};<2.4mm,5.7mm>*{^1}**@{},
    <-0.38mm,-0.2mm>*{};<-2.8mm,2.5mm>*{^2}**@{},
    \end{xy}}\Ea
    \Ea
\right.
\Eeq
Thus a representation,
$$
\rho: \LBcd \lon \cE nd_V
$$
of this prop in a differential graded (dg, for short) vector space $V$ is uniquely determined
by the values of $\rho$ on the generators,
$$
\rho \left(
 \begin{xy}
 <0mm,-0.55mm>*{};<0mm,-2.5mm>*{}**@{-},
 <0.5mm,0.5mm>*{};<2.2mm,2.2mm>*{}**@{-},
 <-0.48mm,0.48mm>*{};<-2.2mm,2.2mm>*{}**@{-},
 <0mm,0mm>*{\circ};<0mm,0mm>*{}**@{},
 \end{xy}\right): V[-c]\rar \odot^2(V[-c])[1], \  \ \ \
\left(
 \begin{xy}
 <0mm,0.66mm>*{};<0mm,3mm>*{}**@{-},
 <0.39mm,-0.39mm>*{};<2.2mm,-2.2mm>*{}**@{-},
 <-0.35mm,-0.35mm>*{};<-2.2mm,-2.2mm>*{}**@{-},
 <0mm,0mm>*{\circ};<0mm,0mm>*{}**@{},
 \end{xy}
 \right): \odot^2(V[d]) \rar V[1+d],
$$
which equip  $V$ with  (degree shifted) dg Lie algebra and Lie coalgebra structures satisfying the Drinfeld compatibility condition (which is assured by the vanishing under $\rho$ of the bottom graph in $\cR$).

\sip

The minimal resolution of the prop $\LBcd$
is a free cofibrant prop  $\HoLBcd$ generated by the $\bS$-bimodule ${E}=\{{E}(m,n)\}$ with ${E}(m,n)\neq 0$ only for $m+n\geq 3$ and $m,n\geq 1$,
\Beq\label{2: symmetries of HoLiebcd corollas}
{E}(m,n):=\sgn_m^{\ot |c|}\ot \sgn_n^{|d|}[c[m-1] + d[n-1]-1]=:\text{span}\left\langle
\Ba{c}\resizebox{17mm}{!}{\begin{xy}
 <0mm,0mm>*{\circ};<0mm,0mm>*{}**@{},
 <-0.6mm,0.44mm>*{};<-8mm,5mm>*{}**@{-},
 <-0.4mm,0.7mm>*{};<-4.5mm,5mm>*{}**@{-},
 <0mm,0mm>*{};<1mm,5mm>*{\ldots}**@{},
 <0.4mm,0.7mm>*{};<4.5mm,5mm>*{}**@{-},
 <0.6mm,0.44mm>*{};<8mm,5mm>*{}**@{-},
   <0mm,0mm>*{};<-10.5mm,5.9mm>*{^{\sigma(1)}}**@{},
   <0mm,0mm>*{};<-4mm,5.9mm>*{^{\sigma(2)}}**@{},
   <0mm,0mm>*{};<10.0mm,5.9mm>*{^{\sigma(m)}}**@{},
 <-0.6mm,-0.44mm>*{};<-8mm,-5mm>*{}**@{-},
 <-0.4mm,-0.7mm>*{};<-4.5mm,-5mm>*{}**@{-},
 <0mm,0mm>*{};<1mm,-5mm>*{\ldots}**@{},
 <0.4mm,-0.7mm>*{};<4.5mm,-5mm>*{}**@{-},
 <0.6mm,-0.44mm>*{};<8mm,-5mm>*{}**@{-},
   <0mm,0mm>*{};<-10.5mm,-6.9mm>*{^{\tau(1)}}**@{},
   <0mm,0mm>*{};<-4mm,-6.9mm>*{^{\tau(2)}}**@{},
   <0mm,0mm>*{};<10.0mm,-6.9mm>*{^{\tau(n)}}**@{},
 \end{xy}}\Ea
=(-1)^{c|\sigma|+d|\tau|}
\Ba{c}\resizebox{14mm}{!}{\begin{xy}
 <0mm,0mm>*{\circ};<0mm,0mm>*{}**@{},
 <-0.6mm,0.44mm>*{};<-8mm,5mm>*{}**@{-},
 <-0.4mm,0.7mm>*{};<-4.5mm,5mm>*{}**@{-},
 <0mm,0mm>*{};<-1mm,5mm>*{\ldots}**@{},
 <0.4mm,0.7mm>*{};<4.5mm,5mm>*{}**@{-},
 <0.6mm,0.44mm>*{};<8mm,5mm>*{}**@{-},
   <0mm,0mm>*{};<-8.5mm,5.5mm>*{^1}**@{},
   <0mm,0mm>*{};<-5mm,5.5mm>*{^2}**@{},
   <0mm,0mm>*{};<4.5mm,5.5mm>*{^{m\hspace{-0.5mm}-\hspace{-0.5mm}1}}**@{},
   <0mm,0mm>*{};<9.0mm,5.5mm>*{^m}**@{},
 <-0.6mm,-0.44mm>*{};<-8mm,-5mm>*{}**@{-},
 <-0.4mm,-0.7mm>*{};<-4.5mm,-5mm>*{}**@{-},
 <0mm,0mm>*{};<-1mm,-5mm>*{\ldots}**@{},
 <0.4mm,-0.7mm>*{};<4.5mm,-5mm>*{}**@{-},
 <0.6mm,-0.44mm>*{};<8mm,-5mm>*{}**@{-},
   <0mm,0mm>*{};<-8.5mm,-6.9mm>*{^1}**@{},
   <0mm,0mm>*{};<-5mm,-6.9mm>*{^2}**@{},
   <0mm,0mm>*{};<4.5mm,-6.9mm>*{^{n\hspace{-0.5mm}-\hspace{-0.5mm}1}}**@{},
   <0mm,0mm>*{};<9.0mm,-6.9mm>*{^n}**@{},
 \end{xy}}\Ea
 \right\rangle_{ \forall \sigma\in \bS_m, \forall\tau\in \bS_n}
\Eeq
The differential on $\HoLBcd$ is given on the generators by
\Beq\label{LBcd_infty}
\delta
\Ba{c}\resizebox{14mm}{!}{\begin{xy}
 <0mm,0mm>*{\circ};<0mm,0mm>*{}**@{},
 <-0.6mm,0.44mm>*{};<-8mm,5mm>*{}**@{-},
 <-0.4mm,0.7mm>*{};<-4.5mm,5mm>*{}**@{-},
 <0mm,0mm>*{};<-1mm,5mm>*{\ldots}**@{},
 <0.4mm,0.7mm>*{};<4.5mm,5mm>*{}**@{-},
 <0.6mm,0.44mm>*{};<8mm,5mm>*{}**@{-},
   <0mm,0mm>*{};<-8.5mm,5.5mm>*{^1}**@{},
   <0mm,0mm>*{};<-5mm,5.5mm>*{^2}**@{},
   <0mm,0mm>*{};<4.5mm,5.5mm>*{^{m\hspace{-0.5mm}-\hspace{-0.5mm}1}}**@{},
   <0mm,0mm>*{};<9.0mm,5.5mm>*{^m}**@{},
 <-0.6mm,-0.44
 mm>*{};<-8mm,-5mm>*{}**@{-},
 <-0.4mm,-0.7mm>*{};<-4.5mm,-5mm>*{}**@{-},
 <0mm,0mm>*{};<-1mm,-5mm>*{\ldots}**@{},
 <0.4mm,-0.7mm>*{};<4.5mm,-5mm>*{}**@{-},
 <0.6mm,-0.44mm>*{};<8mm,-5mm>*{}**@{-},
   <0mm,0mm>*{};<-8.5mm,-6.9mm>*{^1}**@{},
   <0mm,0mm>*{};<-5mm,-6.9mm>*{^2}**@{},
   <0mm,0mm>*{};<4.5mm,-6.9mm>*{^{n\hspace{-0.5mm}-\hspace{-0.5mm}1}}**@{},
   <0mm,0mm>*{};<9.0mm,-6.9mm>*{^n}**@{},
 \end{xy}}\Ea
\ \ = \ \
 \sum_{[1,\ldots,m]=I_1\sqcup I_2\atop
 {|I_1|\geq 0, |I_2|\geq 1}}
 \sum_{[1,\ldots,n]=J_1\sqcup J_2\atop
 {|J_1|\geq 1, |J_2|\geq 1}
}\hspace{0mm}
\pm
\Ba{c}\resizebox{22mm}{!}{ \begin{xy}
 <0mm,0mm>*{\circ};<0mm,0mm>*{}**@{},
 <-0.6mm,0.44mm>*{};<-8mm,5mm>*{}**@{-},
 <-0.4mm,0.7mm>*{};<-4.5mm,5mm>*{}**@{-},
 <0mm,0mm>*{};<0mm,5mm>*{\ldots}**@{},
 <0.4mm,0.7mm>*{};<4.5mm,5mm>*{}**@{-},
 <0.6mm,0.44mm>*{};<12.4mm,4.8mm>*{}**@{-},
     <0mm,0mm>*{};<-2mm,7mm>*{\overbrace{\ \ \ \ \ \ \ \ \ \ \ \ }}**@{},
     <0mm,0mm>*{};<-2mm,9mm>*{^{I_1}}**@{},
 <-0.6mm,-0.44mm>*{};<-8mm,-5mm>*{}**@{-},
 <-0.4mm,-0.7mm>*{};<-4.5mm,-5mm>*{}**@{-},
 <0mm,0mm>*{};<-1mm,-5mm>*{\ldots}**@{},
 <0.4mm,-0.7mm>*{};<4.5mm,-5mm>*{}**@{-},
 <0.6mm,-0.44mm>*{};<8mm,-5mm>*{}**@{-},
      <0mm,0mm>*{};<0mm,-7mm>*{\underbrace{\ \ \ \ \ \ \ \ \ \ \ \ \ \ \
      }}**@{},
      <0mm,0mm>*{};<0mm,-10.6mm>*{_{J_1}}**@{},
 <13mm,5mm>*{};<13mm,5mm>*{\circ}**@{},
 <12.6mm,5.44mm>*{};<5mm,10mm>*{}**@{-},
 <12.6mm,5.7mm>*{};<8.5mm,10mm>*{}**@{-},
 <13mm,5mm>*{};<13mm,10mm>*{\ldots}**@{},
 <13.4mm,5.7mm>*{};<16.5mm,10mm>*{}**@{-},
 <13.6mm,5.44mm>*{};<20mm,10mm>*{}**@{-},
      <13mm,5mm>*{};<13mm,12mm>*{\overbrace{\ \ \ \ \ \ \ \ \ \ \ \ \ \ }}**@{},
      <13mm,5mm>*{};<13mm,14mm>*{^{I_2}}**@{},
 <12.4mm,4.3mm>*{};<8mm,0mm>*{}**@{-},
 <12.6mm,4.3mm>*{};<12mm,0mm>*{\ldots}**@{},
 <13.4mm,4.5mm>*{};<16.5mm,0mm>*{}**@{-},
 <13.6mm,4.8mm>*{};<20mm,0mm>*{}**@{-},
     <13mm,5mm>*{};<14.3mm,-2mm>*{\underbrace{\ \ \ \ \ \ \ \ \ \ \ }}**@{},
     <13mm,5mm>*{};<14.3mm,-4.5mm>*{_{J_2}}**@{},
 \end{xy}}\Ea
\Eeq
where the signs on the r.h.s\ are uniquely fixed for $c+d\in 2\Z$ by the fact that they all equal to $+1$ if $ c$ and $d$ are even integers, and for $c+d\in 2\Z+1$ the signs are given explicitly in
\cite{Me1}. Note that the props $\HoLBcd$ and $\HoLB_{d,c}$ are canonically isomorphic to each other via the flow reversing on the generating graphs.

\sip

A representation of $\HoLBcd$ in a finite-dimensional vector space $V$ can be identified with
a degree $c+d+1$ element $\pi$ in the completed graded commutative algebra
$$
\pi=\sum_{m,n\geq 1\atop m+n\geq 3}=\prod_{m,n\geq 1,m+n\geq 3}
\Hom(\odot^n (V[d]), \odot^m (V[-c])\subset \prod_{k\geq 0}\odot^k (V^*[-d]\oplus V[-c])
$$
equipped with the obvious Poisson type Lie bracket of degree $-c-d$.

\subsection{Non-cofibrant extensions of $\HoLBcd$}
Consider a dg prop  $\sHoLBcd$
generated by the $\bS$-bimodule $E^{\star}=\{E^\star(m,n)\}_{m,n\geq 0}$ with {\em all}\,  $E^\star(m,n)$ non-zero and given by the same formula as in (\ref{2: symmetries of HoLiebcd corollas}).
The
 differential $\delta^{\star}$ on $\suHoLBcd$ is given formally by the formula
(\ref{LBcd_infty}) with the summation over partitions of the sets $[m]$ and $[n]$ appropriately extended,
\Beq\label{2: differential in HoLBcd^star}
\delta^{\star}
\Ba{c}\resizebox{14mm}{!}{\begin{xy}
 <0mm,0mm>*{\circ};<0mm,0mm>*{}**@{},
 <-0.6mm,0.44mm>*{};<-8mm,5mm>*{}**@{-},
 <-0.4mm,0.7mm>*{};<-4.5mm,5mm>*{}**@{-},
 <0mm,0mm>*{};<-1mm,5mm>*{\ldots}**@{},
 <0.4mm,0.7mm>*{};<4.5mm,5mm>*{}**@{-},
 <0.6mm,0.44mm>*{};<8mm,5mm>*{}**@{-},
   <0mm,0mm>*{};<-8.5mm,5.5mm>*{^1}**@{},
   <0mm,0mm>*{};<-5mm,5.5mm>*{^2}**@{},
   <0mm,0mm>*{};<4.5mm,5.5mm>*{^{m\hspace{-0.5mm}-\hspace{-0.5mm}1}}**@{},
   <0mm,0mm>*{};<9.0mm,5.5mm>*{^m}**@{},
 <-0.6mm,-0.44mm>*{};<-8mm,-5mm>*{}**@{-},
 <-0.4mm,-0.7mm>*{};<-4.5mm,-5mm>*{}**@{-},
 <0mm,0mm>*{};<-1mm,-5mm>*{\ldots}**@{},
 <0.4mm,-0.7mm>*{};<4.5mm,-5mm>*{}**@{-},
 <0.6mm,-0.44mm>*{};<8mm,-5mm>*{}**@{-},
   <0mm,0mm>*{};<-8.5mm,-6.9mm>*{^1}**@{},
   <0mm,0mm>*{};<-5mm,-6.9mm>*{^2}**@{},
   <0mm,0mm>*{};<4.5mm,-6.9mm>*{^{n\hspace{-0.5mm}-\hspace{-0.5mm}1}}**@{},
   <0mm,0mm>*{};<9.0mm,-6.9mm>*{^n}**@{},
 \end{xy}}\Ea
\ \ = \ \
 \sum_{[1,\ldots,m]=I_1\sqcup I_2\atop
 {|I_1|\geq 0, |I_2|\geq 0}}
 \sum_{[1,\ldots,n]=J_1\sqcup J_2\atop
 {|J_1|\geq 0, |J_2|\geq 0}
}\hspace{0mm}
\pm
\Ba{c}\resizebox{22mm}{!}{ \begin{xy}
 <0mm,0mm>*{\circ};<0mm,0mm>*{}**@{},
 <-0.6mm,0.44mm>*{};<-8mm,5mm>*{}**@{-},
 <-0.4mm,0.7mm>*{};<-4.5mm,5mm>*{}**@{-},
 <0mm,0mm>*{};<0mm,5mm>*{\ldots}**@{},
 <0.4mm,0.7mm>*{};<4.5mm,5mm>*{}**@{-},
 <0.6mm,0.44mm>*{};<12.4mm,4.8mm>*{}**@{-},
     <0mm,0mm>*{};<-2mm,7mm>*{\overbrace{\ \ \ \ \ \ \ \ \ \ \ \ }}**@{},
     <0mm,0mm>*{};<-2mm,9mm>*{^{I_1}}**@{},
 <-0.6mm,-0.44mm>*{};<-8mm,-5mm>*{}**@{-},
 <-0.4mm,-0.7mm>*{};<-4.5mm,-5mm>*{}**@{-},
 <0mm,0mm>*{};<-1mm,-5mm>*{\ldots}**@{},
 <0.4mm,-0.7mm>*{};<4.5mm,-5mm>*{}**@{-},
 <0.6mm,-0.44mm>*{};<8mm,-5mm>*{}**@{-},
      <0mm,0mm>*{};<0mm,-7mm>*{\underbrace{\ \ \ \ \ \ \ \ \ \ \ \ \ \ \
      }}**@{},
      <0mm,0mm>*{};<0mm,-10.6mm>*{_{J_1}}**@{},
 <13mm,5mm>*{};<13mm,5mm>*{\circ}**@{},
 <12.6mm,5.44mm>*{};<5mm,10mm>*{}**@{-},
 <12.6mm,5.7mm>*{};<8.5mm,10mm>*{}**@{-},
 <13mm,5mm>*{};<13mm,10mm>*{\ldots}**@{},
 <13.4mm,5.7mm>*{};<16.5mm,10mm>*{}**@{-},
 <13.6mm,5.44mm>*{};<20mm,10mm>*{}**@{-},
      <13mm,5mm>*{};<13mm,12mm>*{\overbrace{\ \ \ \ \ \ \ \ \ \ \ \ \ \ }}**@{},
      <13mm,5mm>*{};<13mm,14mm>*{^{I_2}}**@{},
 <12.4mm,4.3mm>*{};<8mm,0mm>*{}**@{-},
 <12.6mm,4.3mm>*{};<12mm,0mm>*{\ldots}**@{},
 <13.4mm,4.5mm>*{};<16.5mm,0mm>*{}**@{-},
 <13.6mm,4.8mm>*{};<20mm,0mm>*{}**@{-},
     <13mm,5mm>*{};<14.3mm,-2mm>*{\underbrace{\ \ \ \ \ \ \ \ \ \ \ }}**@{},
     <13mm,5mm>*{};<14.3mm,-4.5mm>*{_{J_2}}**@{},
 \end{xy}}\Ea
\Eeq
Its ideal $I_0$ generated by all $(m,n)$-corollas\footnote{We often call corollas of type $(0,n)$ (resp.\, $(m,0)$) {\em sources} (resp., {\em targets}). Note that the $(0,0)$ corolla $\bu$ is the unique generator which is both a source and a target. The $(1,1)$-corolla  is often called a {\em passing vertex}.}  with $m=0$ or $n=0$
is differential, and the quotient properad $\HoLBcd^{\star}/I^{0}$ is denoted by
$\HoLBcd^+$ (there exists a general ``plus" endofunctor , $\cP\rar \cP^+$, in the category of props, and the non-cofibrant prop $\HoLBcd^+$ can be understood as the application of that construction to $\HoLBcd$).

 \sip

 The dg prop $\HoLBcd^+$ contains in turn the differential ideal $I^+$ generated by
the $(1,1)$-corolla, and the quotient properad is precisely $\HoLBcd$.

\subsection{Wheeled closures} We refer to \cite{Me1,MMS} for the full details of the wheelification functor, but as we work in this paper only with free props $\cP$ generated by certain $(m,n)$ corollas, $m,n\in \N$, it is very easy to explain what is the wheeled closure
$\cP^\circlearrowright$ of $\cP$: if elements of $\cP$ are obtained in general by glueing output legs of generating corollas to input legs of  other corollas in such a way that directed paths of edges in the resulting directed graph never form a cycle (a ``wheel"), elements of $\cP^\circlearrowright$ are constructed in the same but with with latter restriction dropped. For example,
$$
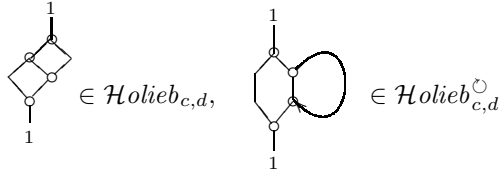
\begin{xy}
 <0mm,-1.3mm>*{};<0mm,-3.5mm>*{}**@{-},
 <0.4mm,0.0mm>*{};<2.4mm,2.1mm>*{}**@{-},
 <-0.38mm,-0.2mm>*{};<-2.8mm,2.5mm>*{}**@{-},
<0mm,-0.8mm>*{\circ};<0mm,0.8mm>*{}**@{},
 <2.96mm,2.4mm>*{\circ};<2.45mm,2.35mm>*{}**@{},
 <2.4mm,2.8mm>*{};<0mm,5mm>*{}**@{-},
 <3.35mm,2.9mm>*{};<5.5mm,5mm>*{}**@{-},
     <0mm,-1.3mm>*{};<0mm,-5.4mm>*{_1}**@{},
%
<-2.8mm,2.5mm>*{};<0mm,5mm>*{\circ}**@{},
<-2.8mm,2.5mm>*{};<0mm,5mm>*{}**@{-},
<2.96mm,5mm>*{};<2.96mm,7.5mm>*{\circ}**@{},
<0.2mm,5.1mm>*{};<2.8mm,7.5mm>*{}**@{-},
<0.2mm,5.1mm>*{};<2.8mm,7.5mm>*{}**@{-},
<5.5mm,5mm>*{};<2.8mm,7.5mm>*{}**@{-},
<2.9mm,7.5mm>*{};<2.9mm,10.5mm>*{}**@{-},
<2.9mm,7.5mm>*{};<2.9mm,11.1mm>*{^1}**@{},
    \end{xy}\in \HoLBcd,\ \ \ \
    \begin{xy}
 <0mm,2.47mm>*{};<0mm,-0.5mm>*{}**@{-},
 <0.5mm,3.5mm>*{};<2.2mm,5.2mm>*{}**@{-},
 <-0.48mm,3.48mm>*{};<-2.2mm,5.2mm>*{}**@{-},
 <0mm,3mm>*{\circ};<0mm,3mm>*{}**@{},
  <0mm,-0.8mm>*{\circ};<0mm,-0.8mm>*{}**@{},
<0mm,-0.8mm>*{};<-2.2mm,-3.5mm>*{}**@{-},
   <0mm,0mm>*{};<-2.5mm,10.1mm>*{^1}**@{},
 <-2.5mm,5.7mm>*{\circ};<0mm,0mm>*{}**@{},
<-2.5mm,5.7mm>*{};<-2.5mm,9.4mm>*{}**@{-},
<-2.5mm,5.7mm>*{};<-5mm,3mm>*{}**@{-},
<-5mm,3mm>*{};<-5mm,-0.8mm>*{}**@{-},
 <-2.5mm,-4.2mm>*{\circ};<0mm,3mm>*{}**@{},
 <-2.8mm,-3.6mm>*{};<-5mm,-0.8mm>*{}**@{-},
 <-2.5mm,-4.6mm>*{};<-2.5mm,-7.3mm>*{}**@{-},
  <0mm,0mm>*{};<-2.5mm,-8.9mm>*{_1}**@{},
 (0.4,3.6)*{}
   \ar@{->}@(ur,dr) (0.1,-0.6)*{}
\end{xy}\in \HoLBcd^\circlearrowright
$$
where the orientation on edges is assumed to flow from bottom to the top unless explicitly shown.
Clearly, $\cP$ is a subprop of its wheeled closure $\cP^\circlearrowright$.
It makes sense to talk about representations of ordinary props  in any vector space (finite- or infinite-dimensional), while their wheeled closures can be represented, in general, only in {\em finite-dimensional}\, vector spaces $V$ as graphs with wheels induce trace operations of the form
$V\ot \Hom(V,\K) \rar \K$ which have no sense in infinite dimensions.

\sip

The wheeled closures of $\sHoLBcd$ and $\HoLBcd^+$ are denoted by $\swHoLBcd$ and $\HoLBcd^{+\circlearrowright}$ respectively.

\sip

 Denote by $\widehat{\HoLB}_{c,d}^{\atop \bigstar\circlearrowright}$ (resp., $\widehat{\HoLB}_{c,d}^{\atop +\circlearrowright}$) the vertex completion of the prop  $\swHoLBcd$ (resp., $\widehat{\HoLB}_{c,d}^{\atop +\circlearrowright}$). One must be careful about definitions of representations of these completed props, but for our purposes the following remark will be enough: given any representation of the prop $\swHoLB_{1,0}$ in a finite-dimensional dg vector space $V$,
 $$
 \rho: \swHoLBcd \lon \cE nd_V
 $$
 that is, a formal Poisson structure $\pi\in \cT_{poly}\cM$ on $V^*$ viewed as a formal graded manifold, there is an associated {\em continuous}\, morphism of the topological props
 $$
 \hat{\rho}: \widehat{\HoLB}_{1,0}^{\atop \bigstar\circlearrowright}\lon \cE nd_V[[\hbar]]
 $$
whose value on any generating corolla $e$ of $\widehat{\HoLB}_{1,0}$ is equal to
$\hbar \rho(e)$, that is, a formal Poisson structure $\hbar \pi\in \cT_{poly}\cM[[\hbar]]$. Here $\hbar$ is any formal parameter of homological
degree zero (``Planck constant").

\sip

\subsubsection{{\bf Proposition}} \label{2: prop on acyclicity of HoLBs+,+} {\em (i) The dg subprop of $\widehat{\HoLB}_{c,d}^{\atop +\circlearrowright}$
spanned by graphs with at least one ingoing or at least one outgoing
legs is acyclic while its complement $\widehat{\HoLB}_{c,d}^{\atop +\circlearrowright}(0,0)$
has non-trivial cohomology which is equal to $H^\bu(\widehat{\HoLB}_{c,d}^{\atop \circlearrowright}(0,0), \delta)$.

\sip

(ii) The dg prop   $\widehat{\HoLB}_{c,d}^{\atop \bigstar\circlearrowright}$ is acyclic.}

\begin{proof}
Consider a filtration of the complex $(\widehat{\HoLB}_{c,d}^{\atop +\circlearrowright}, \delta^+)$ by the number of vertices of valency $\geq 3$.
The induced differential on the associated graded attaches to each leg the $(1,1)$-corolla. We can consider another filtration such that the induced differential attaches $(1,1)$-corolla only to the input (or output) leg labelled by number 1. This complex is obviously acyclic. This proves the claim for the required subprop of  $\widehat{\HoLB}_{c,d}^{\atop +\circlearrowright}$. The second claim about non-triviality of $H^\bu(\widehat{\HoLB}_{c,d}^{\atop +\circlearrowright}(0,0), \delta^+)$ follows from the direct examples of non-trivial cohomology classes such as (see \cite{Me1})
$$
\Ba{c}\resizebox{18mm}{!}{ \begin{xy}
<-5mm,5mm>*{\bullet};
<-5mm,5mm>*{};<-5mm,8mm>*{}**@{-},
<-5mm,5mm>*{};<-7mm,4mm>*{}**@{-},
 <0mm,0mm>*{\bullet};
<0mm,0mm>*{};<-5mm,-5mm>*{}**@{-},
 <0mm,0mm>*{};<-5mm,5mm>*{}**@{-},
<0mm,0mm>*{};<1.5mm,1.5mm>*{}**@{-},
 <0mm,0mm>*{};<1.5mm,-1.5mm>*{}**@{-};
<-5mm,-5mm>*{\bullet};
<-5mm,-5mm>*{};<-7mm,-2mm>*{}**@{-},
<-5mm,-5mm>*{};<-5mm,-8mm>*{}**@{-},
   \ar@{->}@(ul,dl) (-5.0,8.0)*{};(-7.0,4.0)*{},
   \ar@{->}@(ur,dr) (1.5,1.5)*{};(1.5,-1.5)*{},
   \ar@{->}@(ul,dl) (-7.0,-2.0)*{};(-5.0,-8.0)*{},
\end{xy}}\Ea
 -
\Ba{c}\resizebox{18mm}{!}{  \begin{xy}
<0mm,4mm>*{\bullet};
<0mm,4mm>*{};<0mm,8mm>*{}**@{-},
<0mm,4mm>*{};<3mm,2mm>*{}**@{-},
 <0mm,-5mm>*{\bullet};
<0mm,-5mm>*{};<0mm,4mm>*{}**@{-},
<0mm,-5mm>*{};<2mm,-3mm>*{}**@{-},
<0mm,-5mm>*{};<0mm,-8mm>*{}**@{-},
<-5mm,-5mm>*{};<0mm,4mm>*{}**@{-};
<-5mm,-5mm>*{\bullet};
<-5mm,-5mm>*{};<-7mm,-2mm>*{}**@{-},
<-5mm,-5mm>*{};<-5mm,-8mm>*{}**@{-},
   \ar@{->}@(ur,dr) (0,8.0)*{};(3.0,2.0)*{},
   \ar@{->}@(ur,dr) (2.0,-3.0)*{};(0.0,-8.0)*{},
   \ar@{->}@(ul,dl) (-7.0,-2.0)*{};(-5.0,-8.0)*{},
\end{xy}}\Ea
 +
\Ba{c}\resizebox{18mm}{!}{  \begin{xy}
<0mm,-4mm>*{\bullet};
<0mm,-4mm>*{};<0mm,-8mm>*{}**@{-},
<0mm,-4mm>*{};<3mm,-2mm>*{}**@{-},
 <0mm,5mm>*{\bullet};
<0mm,5mm>*{};<0mm,-4mm>*{}**@{-},
<0mm,5mm>*{};<2mm,3mm>*{}**@{-},
<0mm,5mm>*{};<0mm,8mm>*{}**@{-},
<-5mm,5mm>*{};<0mm,-4mm>*{}**@{-};
<-5mm,5mm>*{\bullet};
<-5mm,5mm>*{};<-7mm,2mm>*{}**@{-},
<-5mm,5mm>*{};<-5mm,8mm>*{}**@{-},
   \ar@{->}@(ur,dr) (3.0,-2.0)*{};(0.0,-8.0)*{},
   \ar@{->}@(ur,dr) (0.0,8.0)*{};(2.0,3.0)*{},
   \ar@{->}@(ul,dl) (-5.0,8.0)*{};(-7.0,2.0)*{},
\end{xy}}\Ea
 \in H(\widehat{\HoLB}_{c,d}^{\atop +\circlearrowright}(0,0), \delta^+)\ \ \forall c,d\in \Z \ \text{with} \ c+d\in 2\Z+1,
$$
or even a simpler one
$$
\Ba{c}\resizebox{20mm}{!}{
 \begin{xy}
<-5mm,5mm>*{\bullet};
<-5mm,5mm>*{};<-5mm,8mm>*{}**@{-},
<-5mm,5mm>*{};<-7mm,4mm>*{}**@{-},
 <0mm,0mm>*{\bullet};
 <0mm,0mm>*{};<-5mm,5mm>*{}**@{-},
<0mm,0mm>*{};<1.5mm,1.5mm>*{}**@{-},
 <0mm,0mm>*{};<1.5mm,-1.5mm>*{}**@{-};
%
   \ar@{->}@(ul,dl) (-5.0,8.0)*{};(-7.0,4.0)*{},
   \ar@{->}@(ur,dr) (1.5,1.5)*{};(1.5,-1.5)*{},
%
\end{xy}}\Ea
 \in H(\widehat{\HoLB}_{c,d}^{\atop +\circlearrowright}(0,0), \delta^+) \ \ \forall c,d\in \Z.
$$
It is easy to see (cf.\ \cite{Wi2}) that graphs containing passing vertices
do not contribute to the cohomology so that
$$
H^\bu(\widehat{\HoLB}_{c,d}^{\atop +\circlearrowright}(0,0), \delta^+)=
H^\bu(\widehat{\HoLB}_{c,d}^{\atop \circlearrowright}(0,0), \delta).
$$

\sip

Consider next a filtration of each complex $(\widehat{\HoLB}_{c,d}^{\atop \bigstar\circlearrowright}(m,n), \delta^{\star})$ with $m+n\geq 1$ by the total number of vertices with no input
edges or no output edges. The induced differential in the associated graded is precisely $\delta^+$ so that the argument as above proves its acyclicity.

\sip

Finally, consider the complex $(\widehat{\HoLB}_{c,d}^{\atop \bigstar\circlearrowright}(0,0), \delta^{\star})$.  Call univalent vertices and passing vertices of graphs from $\widehat{\HoLB}_{c,d}^{\atop \bigstar\circlearrowright}(0,0)$
 {\em stringy}\, ones, and call    maximal connected subgraphs (if any) of a graph $\Ga$ from  $\widehat{\HoLB}_{c,d}^{\atop \bigstar\circlearrowright}(0,0)$  consisting of stringy vertices with at least one vertex univalent  {\em strings}.  The vertices of $\Ga$ which do not belong to strings are called {\em core}\, ones. Thus strings are subgraphs or graphs of the following three types,
\Bi
\item[(i)]
$\Ba{c}
\resizebox{55mm}{!}
{ \xy
(-45,1)*+{_\text{core vertex}}="-1",
(-30,1)*{\bu}="0",
(-20,1)*{\bu}="1",
(-13,1)*{}="2",
(-10,1)*{\ldots},
(-7,1)*{}="3",
(0,1)*{\bu}="4",
(10,1)*{\bu}="5",
%
\ar @{->} "-1";"0" <0pt>
\ar @{->} "0";"1" <0pt>
\ar @{->} "1";"2" <0pt>
\ar @{->} "3";"4" <0pt>
\ar @{->} "4";"5" <0pt>
\endxy}
\Ea$\ \text{$n\geq 0$ stringy vertices (shown as black bullets)}
\item[(ii)]
$\Ba{c}
\resizebox{55mm}{!}
{ \xy
(-45,1)*+{_\text{core vertex}}="-1",
(-30,1)*{\bu}="0",
(-20,1)*{\bu}="1",
(-13,1)*{}="2",
(-10,1)*{\ldots},
(-7,1)*{}="3",
(0,1)*{\bu}="4",
(10,1)*{\bu}="5",
%
\ar @{<-} "-1";"0" <0pt>
\ar @{<-} "0";"1" <0pt>
\ar @{<-} "1";"2" <0pt>
\ar @{<-} "3";"4" <0pt>
\ar @{<-} "4";"5" <0pt>
\endxy}
\Ea$\ \text{$n\geq 0$  stringy vertices}

\item[(iii)]
$\Ba{c}
\resizebox{55mm}{!}
{ \xy
(-40,1)*{\bu}="-1",
(-30,1)*{\bu}="0",
(-20,1)*{\bu}="1",
(-13,1)*{}="2",
(-10,1)*{\ldots},
(-7,1)*{}="3",
(0,1)*+{\bu}="4",
(10,1)*+{\bu}="5",
(20,1)*+{\bu}="6",
\ar @{->} "-1";"0" <0pt>
\ar @{->} "0";"1" <0pt>
\ar @{->} "1";"2" <0pt>
\ar @{->} "3";"4" <0pt>
\ar @{->} "4";"5" <0pt>
\ar @{->} "5";"6" <0pt>
\endxy}
\Ea$ \ \text{$n\geq 1$ stringy vertices}
\Ei
 Consider a (complete, exhaustive, bounded above) filtration of $(\widehat{\HoLB}_{c,d}^{\atop \star\circlearrowright}(0,0), \delta^{\star})$ by the number of core vertices,
$$
F_{-p}\ \text{is generated by graphs with the number of core vertices}\ \geq p.
$$
The differential in the associated graded acts non-trivially  on strings of types (i) and (ii)  (resp., (iii)) with {\em even}\,  (resp., {\em odd}) number of stringy vertices only
by increasing that number by one. Hence the complexes  $C_{(i)}$, $C_{(ii)}$ and  $C_{(iii)}$ generated by strings of type (i), (ii) and (iii) respectively are acyclic.

\sip

If the set of core vertices is empty, we are in the situation of the complex $C_{(iii)}$ so that the associated cohomology vanishes.

\sip

If the set of core vertices is non-empty,  then the associated graded is isomorphic to the unordered tensor product
$$
\bigotimes_{v} \odot^\bu C_{(i)}^v \ot \odot^\bu  C_{(ii)}^v
$$
over the set of core vertices of the graded symmetric tensor algebras of acyclic complexes   $C_{(i)}$ and $C_{(ii)}$, and hence is acyclic itself.

\end{proof}

 \sip

\sip

{
\Large
\section{\bf Deformation complexes of wheeled props and graph complexes}
}

\sip

\subsection{Derivations of wheeled props} A wheeled prop $\cP^\circlearrowright$ in the category of complexes is an $\bS$-bimodule, that is a collection $\{\cP^\circlearrowright(m,n)\}$ of $(\bS_m)^{op}\times \bS_n$ modules, equipped with two basic operations satisfying certain axioms (see \S 2 in \cite{MMS} for full details, or just pictures 5 and 6 in \cite{MMS}):
\Bi
\item[(i)] the horizontal composition (``a map from  disjoint union of two decorated corollas into a single corolla")
    $$
    \Ba{rccc}
    \circ_h: &\cP^\circlearrowright(m_1,n_1) \ot \cP^\circlearrowright(m_2,n_2) &\lon &
    \cP^\circlearrowright(m_1+m_2,n_1+n_2)\\
    & a\ot b & \lon & a\circ_h b
    \Ea
    $$
\item[(ii)] the trace operation defined for any $m,n\geq 1$ and any $i\in [m]$, $j\in [n]$  (``gluing $i$-th output leg to the $j$-in input leg, and then contracting the resulting internal edge"),
    $$
    \Ba{rccc}
    Tr_j^i: &\cP^\circlearrowright(m,n) &\lon &
    \cP^\circlearrowright(m-1,n-1)\\
    & a & \lon & Tr_j^i(a).
    \Ea
    $$
\Ei
The Lie algebra of derivations of $\cP^\circlearrowright$ is defined
 as the vector space $\Der(\cP^\circlearrowright)\hook \Hom_{\bS}(\cP^\circlearrowright,\cP^\circlearrowright)$  of those endomorphisms
 $D: \cP^\circlearrowright \lon \cP^\circlearrowright$ of the $\bS$-bimodule $\cP^\circlearrowright$ which satisfy the two conditions: (i)
 for any  $a,b\in \cP$ one has
$$
D(a\circ_h b) = D(a)\circ_h f(b) + (-1)^{|D||a|} f(a)\circ_h D(b),\ \ \ \
$$
 and (ii) for any $c\in \cP(m,n)$ with $m,n\geq 1$ and any $i\in [m]$ and $j\in [n]$
$$
D(Tr_j^i(c))= Tr_j^i(D(c)).
$$
If the $\delta$ is a differential in the wheeled prop $\cP^\circlearrowright$, then
$\delta$ is a MC element in $\Der(\cP^\circlearrowright)$ so that the latter becomes also a complex with the differential $d=[\delta,\ ]$.

\sip
We are interested in the complex of derivations of the {\em completed}\, (by the number of vertices)
prop $\widehat{\HoLB}_{c,d}^{\atop \bigstar\circlearrowright}$ but abusing notations denote it from now on by $\Der(\swHoLBcd)$ (cf.\ \cite{MW1}). Any derivation of $\widehat{\HoLB}_{c,d}^{\atop \bigstar\circlearrowright}$  is uniquely determined by its values on the generators of the prop $\HoLBcd^{\atop \bigstar \circlearrowright}$. Hence we have isomorphisms of graded vector spaces,
\Beq\label{2: Der(Holieb++) as graphs}
\Der(\swhHoLBcd) =
 \prod_{m,n\geq 0} \left(\widehat{\HoLB}_{c,d}^{\atop \bigstar\circlearrowright}(m,n) \otimes \sgn_m^{\ot |c|}\otimes \sgn_n^{\ot |d|}\right)^{\bS_m\times \bS_m}[1+c(1-m)+d(1-n)].
\Eeq

Thus elements of this complex
can be interpreted as directed (not necessarily connected) graphs which might  have incoming or outgoing legs and wheels, for example
$$
\Ba{c}
\resizebox{15mm}{!}{ \xy
(0,0)*{\bu}="d1",
(10,0)*{\bu}="d2",
(-5,-5)*{}="dl",
(5,-5)*{}="dc",
(15,-5)*{}="dr",
(0,10)*{\bu}="u1",
(10,10)*{\bu}="u2",
(5,15)*{}="uc",
(5,15)*{}="uc",
(15,15)*{}="ur",
(0,15)*{}="ul",
\ar @{<-} "d1";"d2" <0pt>
\ar @{<-} "u1";"d1" <0pt>
\ar @{->} "u1";"u2" <0pt>
\ar @{<-} "u1";"d2" <0pt>
\ar @{->} "u2";"d2" <0pt>
\ar @{<-} "u2";"d1" <0pt>
\endxy}
\Ea
\hspace{-5mm}
\Ba{c}
\resizebox{15mm}{!}{ \xy
(0,0)*{\bu}="d1",
(10,0)*{\bu}="d2",
(-5,-5)*{}="dl",
(5,-5)*{}="dc",
(15,-5)*{}="dr",
(0,10)*{\bu}="u1",
(10,10)*{\bu}="u2",
(5,15)*{}="uc",
(5,15)*{}="uc",
(15,15)*{}="ur",
(0,15)*{}="ul",
\ar @{<-} "d1";"d2" <0pt>
\ar @{<-} "d2";"dc" <0pt>
\ar @{<-} "d2";"dr" <0pt>
\ar @{<-} "u1";"d1" <0pt>
\ar @{->} "u1";"u2" <0pt>
\ar @{<-} "u1";"d2" <0pt>
\ar @{->} "u2";"d2" <0pt>
\ar @{<-} "u2";"d1" <0pt>
\ar @{<-} "uc";"u2" <0pt>
\ar @{<-} "ur";"u2" <0pt>
\ar @{<-} "ul";"u1" <0pt>
\endxy}
\Ea \in \Der(\swHoLB_{c,d})
$$
Its subcomplex spanned by {\em oriented} (i.e.\ with no wheels) directed graphs is precisely the derivation complex of $\Der(\sHoLBcd)$, e.g.
$$
\Ba{c}
\resizebox{15mm}{!}{ \xy
(0,0)*{\bu}="d1",
(10,0)*{\bu}="d2",
(-5,-5)*{}="dl",
(5,-5)*{}="dc",
(15,-5)*{}="dr",
(0,10)*{\bu}="u1",
(10,10)*{\bu}="u2",
(5,15)*{}="uc",
(5,15)*{}="uc",
(15,15)*{}="ur",
(0,15)*{}="ul",
\ar @{<-} "d1";"d2" <0pt>
\ar @{<-} "u1";"d1" <0pt>
\ar @{->} "u1";"u2" <0pt>
\ar @{<-} "u1";"d2" <0pt>
\ar @{<-} "u2";"d2" <0pt>
\ar @{<-} "u2";"d1" <0pt>
\ar @{<-} "ul";"u1" <0pt>
\endxy}
\Ea \in \Der(\sHoLBcd)
$$
Note that the outgoing or ingoing legs (if any) of these graphs
are not assigned particular numerical labels; more precisely, their numerical labels are (skew)symmetrized in accordance with the parity of the integer parameters $c$ and $d$.

\sip

The Lie algebra  $\Der(\swhHoLBcd)$ contains a Maurer-Cartan element
$$
\ga^{\star}:=\sum_{m,n\geq 0}\sum_{[m]=I_1\sqcup I_2, [n]=J_1\sqcup J_2 \atop
 |I_1|, |I_2|, |J_1|, |J_2|\geq 0}
\hspace{-1mm}
\pm
\Ba{c}\resizebox{22mm}{!}{ \begin{xy}
 <0mm,0mm>*{\circ};<0mm,0mm>*{}**@{},
 <-0.6mm,0.44mm>*{};<-8mm,5mm>*{}**@{-},
 <-0.4mm,0.7mm>*{};<-4.5mm,5mm>*{}**@{-},
 <0mm,0mm>*{};<0mm,5mm>*{\ldots}**@{},
 <0.4mm,0.7mm>*{};<4.5mm,5mm>*{}**@{-},
 <0.6mm,0.44mm>*{};<12.4mm,4.8mm>*{}**@{-},
     <0mm,0mm>*{};<-2mm,7mm>*{\overbrace{\ \ \ \ \ \ \ \ \ \ \ \ }}**@{},
     <0mm,0mm>*{};<-2mm,9mm>*{^{I_1}}**@{},
 <-0.6mm,-0.44mm>*{};<-8mm,-5mm>*{}**@{-},
 <-0.4mm,-0.7mm>*{};<-4.5mm,-5mm>*{}**@{-},
 <0mm,0mm>*{};<-1mm,-5mm>*{\ldots}**@{},
 <0.4mm,-0.7mm>*{};<4.5mm,-5mm>*{}**@{-},
 <0.6mm,-0.44mm>*{};<8mm,-5mm>*{}**@{-},
      <0mm,0mm>*{};<0mm,-7mm>*{\underbrace{\ \ \ \ \ \ \ \ \ \ \ \ \ \ \
      }}**@{},
      <0mm,0mm>*{};<0mm,-10.6mm>*{_{J_1}}**@{},
 <13mm,5mm>*{};<13mm,5mm>*{\circ}**@{},
 <12.6mm,5.44mm>*{};<5mm,10mm>*{}**@{-},
 <12.6mm,5.7mm>*{};<8.5mm,10mm>*{}**@{-},
 <13mm,5mm>*{};<13mm,10mm>*{\ldots}**@{},
 <13.4mm,5.7mm>*{};<16.5mm,10mm>*{}**@{-},
 <13.6mm,5.44mm>*{};<20mm,10mm>*{}**@{-},
      <13mm,5mm>*{};<13mm,12mm>*{\overbrace{\ \ \ \ \ \ \ \ \ \ \ \ \ \ }}**@{},
      <13mm,5mm>*{};<13mm,14mm>*{^{I_2}}**@{},
 <12.4mm,4.3mm>*{};<8mm,0mm>*{}**@{-},
 <12.6mm,4.3mm>*{};<12mm,0mm>*{\ldots}**@{},
 <13.4mm,4.5mm>*{};<16.5mm,0mm>*{}**@{-},
 <13.6mm,4.8mm>*{};<20mm,0mm>*{}**@{-},
     <13mm,5mm>*{};<14.3mm,-2mm>*{\underbrace{\ \ \ \ \ \ \ \ \ \ \ }}**@{},
     <13mm,5mm>*{};<14.3mm,-4.5mm>*{_{J_2}}**@{},
 \end{xy}}\Ea,
$$
which corresponds to the differential $\delta^{\star}$ in $\swhHoLBcd$. Hence the differential  in the complex $\Der(\HoLBcd^{\atop \bigstar\circlearrowright})$ is given by
\Beq\label{d in Der(Holieb^+)}
d^{\star} \Ga :=[\ga^{\star},\Ga]=
 \delta^{\star}\Ga
  \pm
  \sum
  \Ba{c}\resizebox{14mm}{!}{\begin{xy}
 <0mm,0mm>*{\bu};<0mm,0mm>*{}**@{},
 <-0.6mm,0.44mm>*{};<-8mm,5mm>*{}**@{-},
 <-0.4mm,0.7mm>*{};<-4.5mm,5mm>*{}**@{-},
 <0mm,0mm>*{};<-1mm,5mm>*{\ldots}**@{},
 <0.4mm,0.7mm>*{};<4.5mm,5mm>*{}**@{-},
 <0.6mm,0.44mm>*{};<10mm,6mm>*{}**@{-},
   <0mm,0mm>*{};<12.0mm,7.5mm>*{\Ga}**@{},
 <-0.6mm,-0.44mm>*{};<-8mm,-5mm>*{}**@{-},
 <-0.4mm,-0.7mm>*{};<-4.5mm,-5mm>*{}**@{-},
 <0mm,0mm>*{};<-1mm,-5mm>*{\ldots}**@{},
 <0.4mm,-0.7mm>*{};<4.5mm,-5mm>*{}**@{-},
 <0.6mm,-0.44mm>*{};<8mm,-5mm>*{}**@{-},
 \end{xy}}\Ea
   \mp
  \sum
  \Ba{c}\resizebox{14mm}{!}{\begin{xy}
 <0mm,0mm>*{\bu};<0mm,0mm>*{}**@{},
 <-0.6mm,0.44mm>*{};<-8mm,5mm>*{}**@{-},
 <-0.4mm,0.7mm>*{};<-4.5mm,5mm>*{}**@{-},
 <0mm,0mm>*{};<-1mm,5mm>*{\ldots}**@{},
 <0.4mm,0.7mm>*{};<4.5mm,5mm>*{}**@{-},
 <0.6mm,0.44mm>*{};<-10mm,-6mm>*{}**@{-},
   <0mm,0mm>*{};<-12.0mm,-7.5mm>*{\Ga}**@{},
 <-0.6mm,-0.44mm>*{};<8mm,5mm>*{}**@{-},
 <-0.4mm,-0.7mm>*{};<-4.5mm,-5mm>*{}**@{-},
 <0mm,0mm>*{};<-1mm,-5mm>*{\ldots}**@{},
 <0.4mm,-0.7mm>*{};<4.5mm,-5mm>*{}**@{-},
 <0.6mm,-0.44mm>*{};<8mm,-5mm>*{}**@{-},
 \end{xy}}\Ea
 \Eeq
 where the differential in the first term,
 $$
 \delta^{\star}\Ga= (-1)^{|\Ga|}\sum_v\Ga\circ_v \xy
 (0,-2)*{\bullet}="a",
(0,2)*{\bu}="b",
\ar @{->} "a";"b" <0pt>
\endxy,
 $$
acts on the vertices of $\Ga$ by formula\footnote{That formula might be understood as a substitution into each vertex $v$ the graph
 $ \xy
 (0,-2)*{\bullet}="a",
(0,2)*{\bu}="b",
\ar @{->} "a";"b" <0pt>
\endxy$ and redistributing all edges of $v$ along the pair of new created vertices in all possible ways.}
 (\ref{d in Der(Holieb^+)}) while in the remaining two terms one attaches  $(m,n+1)$-corollas and, respectively,  $(m+1,n)$-corollas
   to each outgoing leg  (if any), and, respectively each ingoing leg (if any) of $\Ga$, and sums over all $m,n$ satisfying $m,n\geq 0$.

   \sip

   It is often useful (cf.\ \cite{Wi1,MW1}) to include the graph $\uparrow$ without vertices into the complex $\Der(\swHoLBcd)$ and set, in accordance with the above general formula for $d^\star$,
   $$
   d^\star \uparrow \, =
   \sum_{m,n\geq 0}(m-n)
  \overbrace{
  \underbrace{
 \Ba{c}\resizebox{6mm}{!}  {\xy
(0,4.5)*+{...},
(0,-4.5)*+{...},
(0,0)*{\bu}="o",
(-5,5)*{}="1",
(-3,5)*{}="2",
(3,5)*{}="3",
(5,5)*{}="4",
(-3,-5)*{}="5",
(3,-5)*{}="6",
(5,-5)*{}="7",
(-5,-5)*{}="8",
\ar @{-} "o";"1" <0pt>
\ar @{-} "o";"2" <0pt>
\ar @{-} "o";"3" <0pt>
\ar @{-} "o";"4" <0pt>
\ar @{-} "o";"5" <0pt>
\ar @{-} "o";"6" <0pt>
\ar @{-} "o";"7" <0pt>
\ar @{-} "o";"8" <0pt>
\endxy}\Ea
 }_{n\times}
 }^{m\times},
   $$
The derivation $d^\star \uparrow$ corresponds to the universal automorphism of {\em any}\, dg wheeled prop $\cP^\circlearrowright$ which sends every element
$a\in \cP^\circlearrowright(m,n)$ into $\la^{m-n}a$ for any $\la\in \K\setminus 0$.

\sip

It is important to notice that the subspace
$$
\Der(\swHoLBcd)_{conn}\subset \Der(\swHoLBcd)
$$
spanned by {\em connected}\, graphs is a subcomplex\footnote{The meaning of the complex $\Der(\swHoLBcd)_{conn}$ is that it describes derivations of $\widehat{\HoLB}_{c,d}^{\atop \bigstar\circlearrowright}$ as a {\em properad}\, rather than as a prop.}, and that there is an canonical isomorphism of complexes
\Beq\label{3: Der in terms of Der_connected}
\Der(\swHoLBcd)=\left(\widehat{\odot^\bu} \left(\Der(\swHoLBcd)_{conn}[-1-c-d]\right)\right)[1+c+d]
\Eeq
As the (completed) symmetric tensor product functor is exact, it is enough to compute the cohomology of the subcomplex $\Der(\swHoLBcd)_{conn}$. We do it in the next section
in terms of the cohomology of certain M.\ Kontsevich graph complexes \cite{Ko} which are reminded in the next subsection.


\subsection{Reminder on graph complexes} A {\em graph}\, $\Ga$ is a 1-dimensional $CW$ complex whose 0-cells are called {\em vertices}\, and 1-cells are called {\em edges}. The set of vertices of $\Ga$ is denoted by $V(\Ga)$ and the set of edges by $E(\Ga)$. A graph $\Ga$ is called {\em directed}\, if each edge $e\in E(\Ga)$ comes equipped with a fixed orientation. If a vertex $v$ of a directed graph has
$m\geq 0$ outgoing edges and $n\geq 0$ incoming edges, then we say that $v$ is an $(m,n)$-{\em vertex}.
A $(1,1)$-vertex is called {\em passing}.

\sip

Let $\sG_{n,l}$ be the set of directed graphs $\Ga$ with $n$ vertices and $l$  edges such that
some bijections $V(\Ga)\rar [n]$ and $E(\Ga)\rar [l]$ are fixed, i.e.\ every edges and every vertex of $\Ga$ has a fixed numerical label. There is
a natural right action of the group $\bS_n \times  \bS_l$ on the set $\sG_{n,l}$ with $\bS_n$ acting by relabeling the vertices and  $\bS_l$ by relabeling the
edges. 
Consider a graded vector space  (``directed full graph complex")
$$
\mathsf{dFGC}_d= \prod_{l\geq 0}\prod_{n\geq 1} \K \langle \sG_{n,l}\rangle \ot_{\bS_n\times \bS_l} \left(\sgn_n^{\ot |d|} \ot \sgn_l^{|d-1|}\right) [d(1-n) + l(d-1)]
$$
 This space is spanned by directed graph with no numerical labels on vertices and edges but with a choice of an {\em orientation}: for $d$ even (resp., odd)  this is a choice of ordering of edges (resp., vertices) up to an even permutation. This graded vector space has a Lie algebra structure with
$$
[\Ga_1,\Ga_2]:= \sum_{v\in V(\Ga)} \Ga_1\circ_v \Ga_2 - (-1)^{|\Ga_1||\Ga_2|}
\Ga_2\circ_v \Ga_1
$$
where  $\Ga_1\circ_v \Ga_2$ is defined by substituting the graph $\Ga_2$ into the vertex $v$ of $\Ga_1$ and taking a sum over  re-attachments of dangling edges (attached earlier to $v$) to vertices of $\Ga_2$
in all possible ways.
It is easy to see that the degree 1 graph $\xy
 (0,0)*{\bullet}="a",
(5,0)*{\bu}="b",
\ar @{->} "a";"b" <0pt>
\endxy$ in $\mathsf{dFGC}_d$ is a Maurer-Cartan element, so that one can make
the latter into a complex with the differential
 $$
 \delta:= [\xy
 (0,0)*{\bullet}="a",
(5,0)*{\bu}="b",
\ar @{->} "a";"b" <0pt>
\endxy ,\ ].
 $$
 The complex $\mathsf{dFGC}_d$ contains a subcomplex $\mathsf{FGC}^{or}$
 spanned by {\em oriented graphs}, that is, graphs with no closed paths of directed edges (``wheels").

 \sip

 One can define an {\em undirected}\, full  graph complex as
 $$
\mathsf{FGC}_d= \prod_{l\geq 0}\prod_{n\geq 1} \K \langle \sG_{n,l}\rangle \ot_{\bS_n\times  (\bS_l \ltimes
 (\bS_2)^l)} \left(\sgn_n^{\ot |d|} \ot (\sgn_l^{|d-1|}\ot (\sgn_2^{|d|})^{\ot l})\right) [d(1-n) + l(d-1)]
$$
where the group $(\bS_2)^l$ acts on edges by reversing their directions. This graph complex is spanned by graphs with directions on edges forgotten for $d$ even, and fixed up  to the flip and multiplication by $(-1)$ for $d$ odd.

\sip

 These dg Lie algebras contain dg subalgebras $\mathsf{dcGC}_{d}\subset \mathsf{dFGC}_d$, $\mathsf{cGC}^{or}\subset \mathsf{FGC}_d^{or}$
  and $\mathsf{cGC}_d\subset \mathsf{FGC}_d$ spanned by {\em connected}\, graphs which in turn contain dg Lie subalgebras $\dcGC_d^{\geq 2}$, $\GC_d^{or,\geq 2}$
 and, respectively, $\GC_d^{\geq 2}$ spanned by graphs with all vertices
  having valency $\geq 2$. The dg Lie algebras $\dcGC_d^{\geq 2}$ and $\mathsf{GC}_d^{or,\geq 2}$ (resp., $\GC_d^{\geq 2}$) contain in turn dg Lie subalgebras
  $\dGC_d$ and $\GC_d^{or}$ spanned by graphs with no passing vertices,  (resp., $\GC_d$ spanned by graphs with all vertices at least trivalent).
The canonical inclusion maps
$$
\dGC_d \lon \dcGC_d^{\geq 2} \lon \dcGC_d, \ \ \GC_d^{or} \lon \GC_d^{or,\geq 2}\lon \mathsf{cGC}^{or}_d
$$
are all quasi-isomorphisms \cite{Wi1,Wi2}.  There is also a canonical morphism of dg Lie algebras
\Beq\label{3: from GC to dGC}
\mathsf{GC}^{\geq 2}_2\lon \mathsf{dcGC}_2,
\Eeq
which sends a graph with no directions on edges into a sum of graphs with all possible directions on edges; it is also a quasi-isomorphism \cite{Wi1}. It was proven in
 in \cite{Wi1,Wi2} that
 $$
 H^\bu(\GC_{d}^{\geq 2}) =H^\bu(\GC_{d+1}^{or})
 $$
 and that
$$
H^0(\dGC_2)=H^0(\GC_2^{\geq 2})=H^0(\GC_3^{or})=\grt_1,
$$
where
$\grt_1$ is the Lie algebra of the Grothendieck-Teichm\"uller group $GRT_1$. It is easy to see that
$H^0(\GC_2^{or})=0$ and $H^0(\GC_3^{\geq 2})=0$.

\sip

One has canonical monomorphisms  of complexes
 $$
\odot^{\bu} \left(\mathsf{dGC}_d [-d]\right)[d] \rar  \mathsf{dFGC}_d,\ \ \
\odot^{\bu} \left(\mathsf{GC}^{\geq 2}_d [-d]\right)[d]\rar  \mathsf{FGC}_d,
  \ \ \
\odot^{\bu}\left(\mathsf{GC}^{or}_d [-d]\right)[d]\rar \mathsf{FGC}^{or}_d
 $$
which are quasi-isomorphisms. Hence it is enough to study only connected graph complexes.

\sip

It is often useful  \cite{MW1,MW3} to consider slightly extended dg Lie algebras,
\Beq\label{3: extended GC}
\dGC_d\oplus \ \K ,\ \ \ \GC_d^{\geq 2}\oplus \K ,\ \ \ \  \GC_d^{or}\oplus \K
\Eeq
where the summand $\K$ is generated by an additional element $\emptyset$ concentrated in degree zero, ``a graph with no vertices and edges", whose
Lie bracket, $[\emptyset, \Ga]$, with an element $\Ga$ of $\GC_{d}^{\geq 2}$ or $\GC_{d}^{or}$ is defined as the multiplication of $\Ga$  by twice the number of its loops (in particular, $\emptyset$ is a cycle with respect to the differential
$\delta$). In this case the zero-th cohomology groups of the first two of these extended complexes for $d=2$ and, respectively, of the last complex for $d=3$ are all equal to the Lie algebra $\grt$ of the ``full" Grothendieck-Teichm\"uller group, rather than to its reduced version $\grt_1$.
This very useful fact prompts us to define the {\em full graph complexes}\, of not necessarily connected graphs as the completed graded symmetric tensor algebras
(\ref{1: fGC in terms of GC}) and (\ref{1: fGC^or in terms of GC^0r}).

\bip

{
\Large
\section{\bf Cohomology of the derivation complex of $\swhHoLBcd$}
}

\sip

\subsection{From directed graph complex to the complex of properadic derivations}
Following \cite{MW1} one notices that there is a natural right action of the dg Lie algebra $\mathsf{dcGC}_{c+d+1}$ on the dg wheeled {\em properad}\, $\widehat{\HoLB}_{c,d}^{\atop \star\circlearrowright}$ by properadic derivations, i.e.\ there is a canonical morphism
of dg Lie algebras,
\Beq\label{3: Morhism F from dcGC to Der^++}
\Ba{rccc}
 F^{\star}\colon & \mathsf{dcGC}_{c+d+1} &\to & \Der(\swhHoLBcd)_{conn}\\
         &   \Ga & \to & F(\Ga)
         \Ea
\Eeq
where the derivation $F(\Ga)$ has, by definition, the following values
on the generators of the completed properad  $\wHoLBcd^{\atop\bigstar\circlearrowright}$
\Beq \label{5:derivation Fstar(Ga)}
\left(\Ba{c}\resizebox{12mm}{!}{\begin{xy}
 <0mm,0mm>*{\circ};<0mm,0mm>*{}**@{},
 <-0.6mm,0.44mm>*{};<-8mm,5mm>*{}**@{-},
 <-0.4mm,0.7mm>*{};<-4.5mm,5mm>*{}**@{-},
 <0mm,0mm>*{};<-1mm,5mm>*{\ldots}**@{},
 <0.4mm,0.7mm>*{};<4.5mm,5mm>*{}**@{-},
 <0.6mm,0.44mm>*{};<8mm,5mm>*{}**@{-},
   <0mm,0mm>*{};<-8.5mm,5.5mm>*{^1}**@{},
   <0mm,0mm>*{};<-5mm,5.5mm>*{^2}**@{},
   <0mm,0mm>*{};<9.0mm,5.5mm>*{^m}**@{},
 <-0.6mm,-0.44mm>*{};<-8mm,-5mm>*{}**@{-},
 <-0.4mm,-0.7mm>*{};<-4.5mm,-5mm>*{}**@{-},
 <0mm,0mm>*{};<-1mm,-5mm>*{\ldots}**@{},
 <0.4mm,-0.7mm>*{};<4.5mm,-5mm>*{}**@{-},
 <0.6mm,-0.44mm>*{};<8mm,-5mm>*{}**@{-},
   <0mm,0mm>*{};<-8.5mm,-6.9mm>*{^1}**@{},
   <0mm,0mm>*{};<-5mm,-6.9mm>*{^2}**@{},
   <0mm,0mm>*{};<9.0mm,-6.9mm>*{^n}**@{},
 \end{xy}}\Ea\right)\cdot F^{\star}(\Ga)
=
 \sum_{s:[n]\rar V(\Ga)\atop \hat{s}:[m]\rar V(\Ga)}  \Ba{c}\resizebox{9mm}{!}  {\xy
 (-6,7)*{^1},
(-3,7)*{^2},
(2.5,7)*{},
(7,7)*{^m},
(-3,-8)*{_2},
(3,-6)*{},
(7,-8)*{_n},
(-6,-8)*{_1},
(0,4.5)*+{...},
(0,-4.5)*+{...},
(0,0)*+{\Ga}="o",
(-6,6)*{}="1",
(-3,6)*{}="2",
(3,6)*{}="3",
(6,6)*{}="4",
(-3,-6)*{}="5",
(3,-6)*{}="6",
(6,-6)*{}="7",
(-6,-6)*{}="8",
\ar @{-} "o";"1" <0pt>
\ar @{-} "o";"2" <0pt>
\ar @{-} "o";"3" <0pt>
\ar @{-} "o";"4" <0pt>
\ar @{-} "o";"5" <0pt>
\ar @{-} "o";"6" <0pt>
\ar @{-} "o";"7" <0pt>
\ar @{-} "o";"8" <0pt>
\endxy}\Ea\ \ \ \ \ \ \forall\ m,n\geq 0,
\Eeq
 with the sum being taken over all ways of attaching the incoming and outgoing legs to the graph $\Ga$.
The image
$$
\delta^{\star}:=F^{\star}\left(\xy
 (0,-2)*{\bullet}="a",
(0,2)*{\bu}="b",
\ar @{->} "a";"b" <0pt>
\endxy\right)
$$
gives us the standard differential (\ref{2: differential in HoLBcd^star}) in $\HoLB_{c,d}^{\star\circlearrowright}$.
The monomorphism $\dGC_{c+d+1}\hook \dcGC_{c+d+1}$ is a quasi-isomorphism so that, from the cohomological viewpoint, it is enough to study the restriction of the above map to the dg Lie subalgebra $\dGC_{c+d+1}$ (we denote this restriction by the same symbol).

  \subsubsection{\bf Theorem}\label{5: Theorem on qis to Def(Holieb star_wheeled)}
{\em For any $c,d\in \Z$ the morphism of dg Lie algebras
\Beq\label{5: F^star from dGC to Der HoLB^star}
F^\star : \dGC_{c+d+1} \lon \Der(\swHoLBcd)_{conn}
\Eeq
is a quasi-isomorphism up to one rescaling class represented by the series
 $$
 r^\star=
  \sum_{m,n\geq 0}(m+n-2)
  \overbrace{
  \underbrace{
 \Ba{c}\resizebox{6mm}{!}  {\xy
(0,4.5)*+{...},
(0,-4.5)*+{...},
(0,0)*{\bu}="o",
(-5,5)*{}="1",
(-3,5)*{}="2",
(3,5)*{}="3",
(5,5)*{}="4",
(-3,-5)*{}="5",
(3,-5)*{}="6",
(5,-5)*{}="7",
(-5,-5)*{}="8",
\ar @{-} "o";"1" <0pt>
\ar @{-} "o";"2" <0pt>
\ar @{-} "o";"3" <0pt>
\ar @{-} "o";"4" <0pt>
\ar @{-} "o";"5" <0pt>
\ar @{-} "o";"6" <0pt>
\ar @{-} "o";"7" <0pt>
\ar @{-} "o";"8" <0pt>
\endxy}\Ea
 }_{n\times}
 }^{m\times}.
$$
}

\begin{proof}

For a graph $\Ga$ in  $\Der(\swhHoLBcd)$
let $V^{\lhd 2}(\Ga)\subset V(\Ga)$ be the subset of univalent vertices and passing vertices, and let $V^{\unrhd 2}(\Ga)$ be its complement, i.e.\ the subset of non-passing vertices of valency $\geq 2$ of $\Ga$.
Consider  the following filtration of the complex
$\Der(\swhHoLBcd)$,
\Beq\label{5: filtration od Der by number of geq 2 nonpass vertices}
F_{-p}(\Der(\swhHoLBcd)):=\text{linear span of graphs $\Ga$ with}\  \# V^{\unrhd 2}(\Ga) \geq p.
\Eeq
For a graph $\Ga\in \dGC_{c+d+1}$ one has $V(\Ga)= V^{\unrhd 2}(\Ga)$ so that an analogous filtration of the l.h.s.\ in  (\ref{3: Morhism F from dcGC to Der^++})  takes the form
\Beq\label{5: filtration ofdGC by number of vertices}
F_{-p}(\Der(\dGC_{c+d+1})):=\text{linear span of graphs $\Ga$ with}\  \#V(\Ga) \geq p.
\Eeq
The morphism (\ref{3: Morhism F from dcGC to Der^++}) respects these filtrations and hence
induces the morphism  of the associated graded complexes (all denoted by the same letters),
\Beq\label{5: gr F}
 \cF^\star: (\mathsf{dcGC}_{c+d+1}, 0) \to  (\Der(\swhHoLBcd),\hat{d})
\Eeq
where the induced differential in the l.h.s.\ is trivial
while the induced differential in the r.h.s.\ is given by
 $$
 \hat{d}\Ga=\hat{\delta}\Ga \ \  \pm
  \sum_{\text{in-legs}\atop \text{of}\ \Ga}\left(
  \Ba{c}\resizebox{1.5mm}{!}{\begin{xy}
 <0mm,0mm>*{\bu};
 <0.0mm,7.5mm>*{\Ga};
 <0.0mm,0.44mm>*{};<0mm,5.5mm>*{}**@{-},
 \end{xy}}\Ea
 \pm
 \Ba{c}\resizebox{1.5mm}{!}{\begin{xy}
 <0mm,0mm>*{\bu};
 <0.0mm,7.5mm>*{\Ga};
 <0.0mm,0.44mm>*{};<0mm,5.5mm>*{}**@{-},
 <0.0mm,0.0mm>*{};<0mm,-5.5mm>*{}**@{-},
 \end{xy}}\Ea
 \right)
   \pm\sum_{\text{out-legs}\atop \text{of}\ \Ga}\left(
   \Ba{c}\resizebox{1.5mm}{!}{\begin{xy}
 <0mm,0mm>*{\bu};
 <0.0mm,-7.5mm>*{\Ga};
 <0.0mm,-0.44mm>*{};<0mm,-5.5mm>*{}**@{-},
 \end{xy}}\Ea
 \pm
 \Ba{c}\resizebox{1.5mm}{!}{\begin{xy}
 <0mm,0mm>*{\bu};
 <0.0mm,-7.5mm>*{\Ga};
 <0.0mm,-0.44mm>*{};<0mm,-5.5mm>*{}**@{-},
 <0.0mm,0.0mm>*{};<0mm,5.5mm>*{}**@{-},
 \end{xy}}\Ea
 \right)
  $$
where
$$
\hat{\delta}\Ga=\delta^\star \Ga \bmod \text{terms  creating new $(m,n)$-corollas with $m\geq 2$ or $n\geq 2$}.
$$

Let us call ingoing or outgoing legs (if any) of graphs from (\ref{2: Der(Holieb++) as graphs}) {\em hairs}\, and
consider the following complete, exhaustive and bounded above filtration of  both sides of the arrow in
 (\ref{5: gr F})
\Beq\label{5: filtration by hairs + s + t}
F'_{-p}(\Der(\swhHoLBcd)):=\text{span of graphs with  $\#$hairs + $\#$univalent sources + $\#$univalent targets}\ \geq p.
\Eeq
and
$$
F'_{-p}(\dGC_{c+d+1}):=\left\{\Ba{cl} \dGC_{c+d+1} & \text{for}\ p\leq 0\\
0  & \text{for}\ p\geq 1
\Ea
\right.
$$
 Note that the unique graph $\bu$ consisting of the zero valent vertex
 is counted twice --- once as a source and once as a target --- so that
$\bu$ belongs to $F'_{-2}(\Der(\swhHoLBcd))$; similarly, the derivation $\uparrow$ (see \S 3.1)
is assumed by definition to have two hairs and hence also belongs to $F'_{-2}(\Der(\swhHoLBcd))$.
The map $\cF^\star$ respects both filtrations and hence induces a morphism (denoted by the same letter again) of the associated graded complexes
\Beq\label{5: gr F^st}
 \cF^\star: (\mathsf{dcGC}_{c+d+1}, 0) \to  (\Der(\swhHoLBcd), d_0)=:(C,d_{0})
\Eeq
where the induced differential $d_0$ is given on two exceptional graphs by
 $$
 d_0\, \bu=\Ba{c}\resizebox{1.4mm}{!}  {\xy
 (0,-2.5)*{\bu}="a",
(0,2.5)*{\bu}="b",
\ar @{->} "a";"b" <0pt>
\endxy}\Ea           \ , \ \ \  d_0\uparrow \, = \Ba{c}\resizebox{2.3mm}{!}  {
\xy
 (0,-2.5)*{\bu}="a",
(0,2.5)*{}="b",
\ar @{->} "a";"b" <0pt>
\endxy}\Ea
-
 \Ba{c}\resizebox{1.44mm}{!}  {\xy
 (0,-2.5)*{}="a",
(0,2.5)*{\bu}="b",
\ar @{->} "a";"b" <0pt>
\endxy}\Ea
 $$
 and on all other graphs by the formula
 \Beq\label{5: differential d_0}
 {d}_0\Ga={\delta}^+\Ga \ \  \pm
  \sum_{\text{in-legs}\atop \text{of}\ \Ga}\left(
  \Ba{c}\resizebox{1.5mm}{!}{\begin{xy}
 <0mm,0mm>*{\bu};
 <0.0mm,7.5mm>*{\Ga};
 <0.0mm,0.44mm>*{};<0mm,5.5mm>*{}**@{-},
 \end{xy}}\Ea
 \pm
 \Ba{c}\resizebox{1.5mm}{!}{\begin{xy}
 <0mm,0mm>*{\bu};
 <0.0mm,7.5mm>*{\Ga};
 <0.0mm,0.44mm>*{};<0mm,5.5mm>*{}**@{-},
 <0.0mm,0.0mm>*{};<0mm,-5.5mm>*{}**@{-},
 \end{xy}}\Ea
 \right)
   \pm\sum_{\text{out-legs}\atop \text{of}\ \Ga}\left(
   \Ba{c}\resizebox{1.5mm}{!}{\begin{xy}
 <0mm,0mm>*{\bu};
 <0.0mm,-7.5mm>*{\Ga};
 <0.0mm,-0.44mm>*{};<0mm,-5.5mm>*{}**@{-},
 \end{xy}}\Ea
 \pm
 \Ba{c}\resizebox{1.5mm}{!}{\begin{xy}
 <0mm,0mm>*{\bu};
 <0.0mm,-7.5mm>*{\Ga};
 <0.0mm,-0.44mm>*{};<0mm,-5.5mm>*{}**@{-},
 <0.0mm,0.0mm>*{};<0mm,5.5mm>*{}**@{-},
 \end{xy}}\Ea
 \right)
  \Eeq
  where
\Beqrn
\delta^+\Ga &:=&
\hat{\delta}\Ga \ \text{modulo term creating new univalent vertices}.
\Eeqrn
 By analogy to the proof of Proposition {\ref{2: prop on acyclicity of HoLBs+,+}}, let us call the univalent vertices and passing bivalent vertices
of graphs $\Ga$ from $\Der(\swhHoLBcd)$  {\em stringy}; the   maximal connected subgraphs (if any) of a graph $\Ga$  consisting of stringy vertices with at least one univalent vertex or with at least one hair are called {\em strings}. Let us call the non-passing vertices of valency $\geq 2$ which do not belong to the strings (if any) the {\em core vertices}, and let $\Ga^{core}$ be the full subgraph of $\Ga$ spanned by the core vertices; in principle any graph from the set of generators of $\dcGC_{c+d+1}$ can occur as a core graph $\Ga^{core}$ of some graph $\Ga\in \Der(\swHoLBcd)$.
  A string is a  subgraph (if any) of $\Ga$  of one of the  following eight types (we classify the unique graph in $\Der(\swhHoLBcd)$ consisting consisting of the zero valency vertex $\bu$
  as well as the graph with no vertices $\uparrow$
 as {\em strings}\, as well --- they correspond to the element $\al^\bu_1$  and $\al^{\uparrow\downarrow}_0$ listed below),
\begin{align*}
\al^\bu
_n &\simeq
\Ba{c}
\resizebox{55mm}{!}
{ \xy
(-40,1)*{\bu}="-1",
(-30,1)*{\bu}="0",
(-20,1)*{\bu}="1",
(-13,1)*{}="2",
(-10,1)*{\ldots},
(-7,1)*{}="3",
(0,1)*{\bu}="4",
(10,1)*{\bu}="5",
(20,1)*{\bu}="6",
\ar @{->} "-1";"0" <0pt>
\ar @{->} "0";"1" <0pt>
\ar @{->} "1";"2" <0pt>
\ar @{->} "3";"4" <0pt>
\ar @{->} "4";"5" <0pt>
\ar @{->} "5";"6" <0pt>
\endxy}
\Ea \ \text{$n\geq 1$ stringy vertices}\\
\al^{\uparrow}_n &\simeq
\Ba{c}
\resizebox{55mm}{!}
{ \xy
(-40,1)*{\bu}="-1",
(-30,1)*{\bu}="0",
(-20,1)*{\bu}="1",
(-13,1)*{}="2",
(-10,1)*{\ldots},
(-7,1)*{}="3",
(0,1)*{\bu}="4",
(10,1)*{\bu}="5",
(20,1)*{}="6",
\ar @{->} "-1";"0" <0pt>
\ar @{->} "0";"1" <0pt>
\ar @{->} "1";"2" <0pt>
\ar @{->} "3";"4" <0pt>
\ar @{->} "4";"5" <0pt>
\ar @{->} "5";"6" <0pt>
\endxy}
\Ea \ \text{$n\geq 1$ stringy vertices}\\
\al^{\downarrow}_n & \simeq
\Ba{c}
\resizebox{55mm}{!}
{ \xy
(-40,1)*{}="-1",
(-30,1)*{\bu}="0",
(-20,1)*{\bu}="1",
(-13,1)*{}="2",
(-10,1)*{\ldots},
(-7,1)*{}="3",
(0,1)*{\bu}="4",
(10,1)*{\bu}="5",
(20,1)*{\bu}="6",
\ar @{->} "-1";"0" <0pt>
\ar @{->} "0";"1" <0pt>
\ar @{->} "1";"2" <0pt>
\ar @{->} "3";"4" <0pt>
\ar @{->} "4";"5" <0pt>
\ar @{->} "5";"6" <0pt>
\endxy}
\Ea \ \text{$n\geq 1$ stringy vertices}\\
\al^{\uparrow\downarrow}_n &\simeq
\Ba{c}
\resizebox{55mm}{!}
{ \xy
(-40,1)*{}="-1",
(-30,1)*{\bu}="0",
(-20,1)*{\bu}="1",
(-13,1)*{}="2",
(-10,1)*{\ldots},
(-7,1)*{}="3",
(0,1)*{\bu}="4",
(10,1)*{\bu}="5",
(20,1)*{}="6",
\ar @{->} "-1";"0" <0pt>
\ar @{->} "0";"1" <0pt>
\ar @{->} "1";"2" <0pt>
\ar @{->} "3";"4" <0pt>
\ar @{->} "4";"5" <0pt>
\ar @{->} "5";"6" <0pt>
\endxy}
\Ea \ \text{$n\geq 0$ stringy vertices}
\\
&\\
%
\be^{\bu,\uparrow}_n & \simeq
\Ba{c}
\resizebox{35mm}{!}
{ \xy
(-30,1)*+{_v}*\cir{}="0",
(-20,1)*{\bu}="1",
(-13,1)*{}="2",
(-10,1)*{\ldots},
(-7,1)*{}="3",
(0,1)*{\bu}="4",
(10,1)*{\bu}="5",
\ar @{->} "0";"1" <0pt>
\ar @{->} "1";"2" <0pt>
\ar @{->} "3";"4" <0pt>
\ar @{->} "4";"5" <0pt>
\endxy}
\Ea\ \text{$n\geq 1$ stringy vertices (shown as black bullets)}
\\
\be^{\uparrow}_n & \simeq
\Ba{c}
\resizebox{35mm}{!}
{ \xy
(-30,1)*+{_v}*\cir{}="0",
(-20,1)*{\bu}="1",
(-13,1)*{}="2",
(-10,1)*{\ldots},
(-7,1)*{}="3",
(0,1)*{\bu}="4",
(10,1)*{}="5",
\ar @{->} "0";"1" <0pt>
\ar @{->} "1";"2" <0pt>
\ar @{->} "3";"4" <0pt>
\ar @{->} "4";"5" <0pt>
\endxy}
\Ea\ \text{$n\geq 0$ stringy vertices}
\\
\be^{\bu,\downarrow}_n & \simeq
\Ba{c}
\resizebox{35mm}{!}
{ \xy
(-30,1)*+{_v}*\cir{}="0",
(-20,1)*{\bu}="1",
(-13,1)*{}="2",
(-10,1)*{\ldots},
(-7,1)*{}="3",
(0,1)*{\bu}="4",
(10,1)*{\bu}="5",
\ar @{<-} "0";"1" <0pt>
\ar @{<-} "1";"2" <0pt>
\ar @{<-} "3";"4" <0pt>
\ar @{<-} "4";"5" <0pt>
\endxy}
\Ea\ \text{$n\geq 1$ stringy vertices}
\\
\be^{\downarrow}_n & \simeq
\Ba{c}
\resizebox{35mm}{!}
{ \xy
(-30,1)*+{_v}*\cir{}="0",
(-20,1)*{\bu}="1",
(-13,1)*{}="2",
(-10,1)*{\ldots},
(-7,1)*{}="3",
(0,1)*{\bu}="4",
(10,1)*{}="5",
%
\ar @{<-} "0";"1" <0pt>
\ar @{<-} "1";"2" <0pt>
\ar @{<-} "3";"4" <0pt>
\ar @{<-} "4";"5" <0pt>
\endxy}
\Ea\ \text{$n\geq 0$ stringy vertices}
%
%
%
%
&     \\
& \\
%
\ga_n & \simeq
\Ba{c}
\resizebox{35mm}{!}
{ \xy
(-30,1)*+{_v}*\cir{}="0",
(-20,1)*{\bu}="1",
(-13,1)*{}="2",
(-10,1)*{\ldots},
(-7,1)*{}="3",
(0,1)*{\bu}="4",
(10,1)*+{_w}*\cir{}="5",
\ar @{->} "0";"1" <0pt>
\ar @{->} "1";"2" <0pt>
\ar @{->} "3";"4" <0pt>
\ar @{->} "4";"5" <0pt>
\endxy}
\Ea\ \text{$n\geq 1$ passing vertices (shown as black bullets)}
\end{align*}
where $v$ and $w$  stand for any pair of (not necessary distinct) arbitrary {\em core}\, vertices.
 Note that $\be^{\bu, \uparrow}_0\equiv \be^{\bu, \downarrow}_0$ stand for one and the same element --- a core vertex $v$ with no strings attached.
\sip


 The associated graded complex $C$ in the r.h.s.\ of (\ref{5: gr F^st}) splits into a direct
 sum
 \Beq\label{5: C = C_empty core + C_non-empty core}
 C=C_\text{empty core}\ \ \ \oplus\ \ \  C_\text{non-empty core}
 \Eeq
 where the first (resp., second) summand is spanned by graphs $\Ga$ with the set $V(\Gamma^{core})$ empty (resp., non-empty). Thus
 $$
 C_\text{empty core}:=
 \text{span}\left\langle \al_{n}^{\bu}\ ,
 \al_{n}^{\uparrow}\ ,  \al_{n}^{\downarrow}\ ,  \al_{n}^{\uparrow\downarrow}\ , \uparrow\right\rangle_{n\geq 1}
$$
with the induced differential $d_0$  given on the generators by
$$
d_0 \al_n^{\bu}= \left\{\Ba{ll}
                                   0 & \text{if $n$ is even}\\
                                  \pm \al_{n+1}^{\bu} & \text{if $n$ is odd}   \Ea\right.\ \ \ , \ \ \
                                  d_0 \al_n^{\uparrow\downarrow}= \left\{\Ba{ll} \pm\al_{n+1}^{\uparrow} \pm\al_{n+1}^{\downarrow} \pm\al_{n+1}^{\uparrow\downarrow} & \text{if $n$ is odd}\\
 \pm\al_{n+1}^{\uparrow} \pm\al_{n+1}^{\downarrow} & \text{if $n$ is even}\
                                   \Ea\right.
 $$
 $$
 d_0 \al_n^{\uparrow}= \left\{\Ba{ll} \pm\al_{n+1}^{\bu}  & \text{if $n$ is odd}\\
                                   \pm\al_{n+1}^{\bu} \pm \al_{n+1}^{\uparrow} & \text{if $n$ is even}\
                                   \Ea\right. , \ \ \
    d_0 \al_n^{\downarrow}= \left\{\Ba{ll} \pm\al_{n+1}^{\bu}  & \text{if $n$ is odd}\\
                                   \pm\al_{n+1}^{\bu} \pm \al_{n+1}^{\downarrow} & \text{if $n$ is even}\
                                   \Ea\right.
$$

It is easy to see that the cohomology of this complex is one-dimensional and is equal to the sum of two cycles
\Beq\label{5: graded ass of rescaling class}
\left(\bu +  \Ba{c}\resizebox{2mm}{!}  {
\xy
 (0,-2.5)*{\bu}="a",
(0,2.5)*{}="b",
\ar @{->} "a";"b" <0pt>
\endxy}\Ea\right)
+
\left(\bu +
 \Ba{c}\resizebox{1.24mm}{!}  {\xy
 (0,-2.5)*{}="a",
(0,2.5)*{\bu}="b",
\ar @{->} "a";"b" <0pt>
\endxy}\Ea\right)= 2\bu  + \Ba{c}\resizebox{2mm}{!}  {
\xy
 (0,-2.5)*{\bu}="a",
(0,2.5)*{}="b",
\ar @{->} "a";"b" <0pt>
\endxy}\Ea
+
 \Ba{c}\resizebox{1.24mm}{!}  {\xy
 (0,-2.5)*{}="a",
(0,2.5)*{\bu}="b",
\ar @{->} "a";"b" <0pt>
\endxy}\Ea
\Eeq
(whose difference is a coboundary as $d_0 \uparrow=
\Ba{c}\resizebox{2mm}{!}  {
\xy
 (0,-2.5)*{\bu}="a",
(0,2.5)*{}="b",
\ar @{->} "a";"b" <0pt>
\endxy}\Ea
-
 \Ba{c}\resizebox{1.24mm}{!}  {\xy
 (0,-2.5)*{}="a",
(0,2.5)*{\bu}="b",
\ar @{->} "a";"b" <0pt>
\endxy}\Ea$). This is precisely the representative of the rescaling class $r^\star$.

\sip

Consider next the second complex $(C_\text{non-empty core}, d_0)$.
It decomposes into the completed direct sum (parameterized by arbitrary graphs $\Ga^{core}$ from $\dcGC_{c+d+1}$) of the tensor products of complexes
$$
 C_\text{non-empty core}\simeq \prod_{\Ga^{core}}C_{\Ga^{core}},\ \ \ \  C_{\Gamma^{core}}:= \bigotimes_{v\in V(\Ga^{core})} X_v \bigotimes_{e\in E(\Ga^{core})} X_e
$$
where
\Bi
\item for each edge
$
X_{e}:= \K[0] \oplus \text{span} \left\langle \ga_n\right\rangle_{n\geq 1}
$
with differential given on generators by $d(1\in \K[0])=0$ and $d\ga_n=\pm \ga_{n+1}$ so that $H^\bu(X_e)=\K[0]$; hence the factors $X_e$ can be ignored in the above formula for the complex $C_{\Ga^{core}}$;
\item the complexes $X_v$ can be different for different core vertices $v$ but their classification is rather simple and is discussed next.
\Ei

\sip

For each core vertex $v$ consider two complexes,
$$
C_v^\uparrow:=\text{span}\left\langle \be_{n}^{\bu, \uparrow}\  ,
 \be_{m}^{\uparrow}\ , \right\rangle_{n\geq 1,m\geq 0}, \ \ \ \
C_v^\downarrow:=\text{span}\left\langle \be_{n}^{\bu, \downarrow}\ ,
 \  \be_{m}^{\downarrow}\ \right\rangle_{n\geq 1,m\geq 0}
$$
equipped with the differentials given by

$$
d_v \be_n^{\bu, \uparrow}= \left\{\Ba{ll} \pm \be_{n+1}^{\bu, \uparrow} & \text{if $n$ is even}\\
                                   0 & \text{if $n$ is odd}\
                                   \Ea\right.\ \ , \ \
d_v \be_n^{\uparrow}= \left\{\Ba{ll} \pm\be_{n+1}^{\bu, \uparrow}  & \text{if $n$ is even}\\
                                   \pm\be_{n+1}^{\bu, \uparrow} \pm \be_{n+1}^{\uparrow} & \text{if $n$ is odd}\
                                   \Ea\right.
 $$

 $$
 d_v \be_n^{\bu, \downarrow}= \left\{\Ba{ll} \pm \be_{n+1}^{\bu, \downarrow} & \text{if $n$ is even}\\
                                   0 & \text{if $n$ is odd}\
                                   \Ea\right. , \ \ \
 d_v \be_n^{\downarrow}= \left\{\Ba{ll} \pm\be_{n+1}^{\bu,\downarrow}  & \text{if $n$ is even}\\
                                   \pm\be_{n+1}^{\bu, } \pm \be_{n+1}^{\downarrow} & \text{if $n$ is odd}\
                                   \Ea\right.
 $$

It is easy that that both complexes $C_v^\uparrow$ and $C_v^\downarrow$ are acyclic. Indeed,
consider a filtration of, say, the complex $C_v^\uparrow$ by the number of vertices of the form   $\be^{\bu, \uparrow}_n$; the cohomology of the associated graded complex is 2-dimensional and is spanned by $\be^{\bu, \uparrow}_1$ and $\be^{\uparrow}_0$ with the induced differential given on the generators by the isomorphism $\be^{\uparrow}_0\rar \be^{\bu, \uparrow}_1$.

\mip

Next we have to consider several types of non-empty core graphs.

\sip

{\sc Case 1}: the case $\Ga^{core}=\bu$, the single vertex without any edges. In this case
$$
C_{\Ga^{core}}=\prod_{p+q\geq 3} {\odot}^p C_v^\uparrow \ot
\odot^q C_v^\downarrow \ \ \oplus\ \ \ \odot^2 C_v^\uparrow \ \ \ \oplus \ \ \
\odot^2 C_v^\downarrow.
$$
Due to the acyclicity of the complexes $C_v^\downarrow$ and $C_v^\uparrow$ and exactness
of the (symmetric) tensor product functor, we conclude that $H^\bu(C_{\Ga^{core}})=0$ in this case.

\sip

{\sc Case 2}: the core graph $\Ga^{core}$ contains at least one  vertex $v$ which either has valency one or is a  passing\footnote{Note that a passing vertex in $\Ga^{core}$ can {\em not}\, be passing in $\Ga$.} vertex. Then $C_{\Ga^{core}}$ has the following tensor factor
$$
X_v=\left\{
\Ba{ll}\prod_{p+q\geq 2} \odot^p C_v^\uparrow \ot
\odot^q C_v^\downarrow \ \oplus\ C_v^\uparrow & \text{if}\ |v|=|v|_{out}=1\\
\prod_{p+q\geq 2} \odot^p C_v^\uparrow \ot
\odot^q C_v^\downarrow \ \oplus\ C_v^{\downarrow} & \text{if}\ |v|=|v|_{in}=1\\
\prod_{p+q\geq 1} \odot^p C_v^\uparrow \ot
\odot^q C_v^\downarrow  & \text{if}\ |v|=2, |v|_{in}=|v|_{out}=1
\Ea
\right.
$$
which is in all cases acyclic, $H^\bu(X_v)=0$, so that $H^\bu(C_{\Ga^{core}})=0$.

\sip

Thus we conclude that only generators $\Ga^{core}$ of the subspace $\dGC_{c+d+1}$ can contribute to
$H^\bu(C_\text{non-empty core})$.

\sip

{\sc Case 3}: Consider finally the case when $\Ga^{core}$ is a generator of $\dGC_{c+d+1}$. Then
$$
C_{\Gamma^{core}}:= \bigotimes_{v\in V(\Ga^{core})}\hspace{-2mm}  X_v \bigotimes_{e\in E(\Ga^{core})} \hspace{-2mm} X_e\ \ \
 \text{with}\ \
X_v=\prod_{p,q\geq 0} \odot^p C_v^\uparrow \ot
\odot^q C_v^\downarrow \ \ \ \text{and} \ \ X_{e}:= \K[0] \oplus \text{span} \left\langle \ga_n\right\rangle_{n\geq 1}
$$
We conclude that for each $v\in V(\Ga^{core})$ (resp., each edge $e\in E(\Ga^{core})$) the associated cohomology group $H^\bu(X_v)$ (resp., $H^\bu(X_e)$) is concentrated in degree zero and is equal to $\K$ so that $H^\bu(C_{\Ga^{core}})=\text{span}\left\langle\Ga^{core}\right\rangle$ and hence
$$
H^\bu(C_\text{non-empty core})\simeq \dGC_{c+d+1}.
$$
Moreover, this isomorphism holds true at the level complexes when turning the page  of our spectral sequence.
 By the spectral sequences comparison theorem
 we conclude that
 the map  of the original complexes (\ref{5: F^star from dGC to Der HoLB^star})  is a quasi-isomorphism up to one rescaling class $r^\star$.
\end{proof}

\subsubsection{\bf Remarks}
(i) In terms of representations of the prop $\swHoLB_{0,1}$, that is, in terms of formal Poisson structures, the rescaling class $r^\star$ corresponds to the following universal automorphism
$$
\pi=\sum_{m,n\geq 0} \pi_n^m \lon \pi^{new}=\sum_{m,n\geq 0} \la^{m+n-2} \pi_n^m, \ \ \ \forall \
\la\in \K^\star,
$$
of the set of formal (finite- or infinite-dimensional) Poisson structures.

(ii) In terms of the extended graph complexes (\ref{3: extended GC}) the above result can be re-written as a quasi-isomorphism of dg Lie algebras
$$
\dGC_{c+d+1}\oplus \K \lon \Der(\swHoLBcd)_{conn}
$$
where the generator $\emptyset$ of $\K$ is mapped into the rescaling class. Note
that the l.h.s.\, is {\em not}\, a direct sum of Lie algebras, only of graded vector spaces.

(iii) Composing  the quasi-isomorphism (\ref{3: from GC to dGC}) with the quasi-isomorphism  (\ref{3: Morhism F from dcGC to Der^++}) and using equalities
(\ref{1: fGC in terms of GC}) and (\ref{3: Der in terms of Der_connected}) we obtain a canonical quasi-isomorphism of dg Lie algebras
$$
F^\circlearrowright: \fGC^{\geq 2}_{c+d+1} \lon \Der(\swHoLBcd)
$$
and hence prove the first part of Proposition 1.1.1.

\mip

Similarly one can study the deformation theory of the ordinary (non-wheeled) properad $\sHoLBcd$ and obtain the following result.

\subsubsection{\bf Proposition}\label{5: F from GC^0r to Der(HoLB^star)}
 {\em There is a quasi-isomorphism of dg Lie algebras}
$$
F: \fGC_{c+d+1}^{or}
\lon  \Der(\sHoLBcd)
$$

We skip the details which are identical to the arguments used in the proof of Theorem 4.1.1.


\mip

\bip

{\Large
\section{\bf Classification of universal quantizations \\ of Poisson structures}
}

\mip

\subsection{Polydifferential functor on wheeled props}

There is a polydifferential functor\footnote{In fact in the subsequent paper \cite{MW3} the functor $\f$ was further extended to a {\em polydifferential endofunctor}\, $\caD$ in the category of (wheeled) props such that $\f$ is an operadic part of $\caD$.} \cite{MW2}
$$
\f: \text{\sf Category of dg props} \lon \text{\sf Category of dg operads}
$$
which verbatim extends (on the l.h.s.) to the category of dg wheeled props and has the property
that for any dg (wheeled) prop $\cP$ and its any representation, $\rho: \cP\rar \cE nd_V$,
 in a dg vector space $V$ the associated dg operad $\f(\cP)$ has an associated representation, $\f(\rho): \f(\cP)\rar \cE nd_{\widehat{\odot^\bu} V}$,
 in the (completed) graded commutative algebra $\widehat{\odot^\bu} V$ given in terms of polydifferential (with respect to the standard multiplication in $\widehat{\odot^\bu} V$) operators. We refer to \cite{MW2} for full details and explain here only the explicit structure of the dg operad
 $$
 \f(\widehat{\HoLB}_{0,1}^{\atop \bigstar\circlearrowright})=\left\{
 \f(\widehat{\HoLB}_{0,1}^{\atop \bigstar\circlearrowright})(k)\right\}_{k\geq 0}
 $$
A typical element $a$ in the $\bS_k$-module $\f(\widehat{\HoLB}_{0,1}^{\atop \bigstar\circlearrowright})(k)$ is a linear combination,
$$
\la_1 \hat{e}_1 +  \ldots + \la_N \hat{e}_N, \ \ \ \la_1,\ldots,\la_N\in \K,
$$
where each generator, say,  $\hat{e}_s\in  \f(\widehat{\HoLB}_{0,1}^{\atop \bigstar\circlearrowright})(k)$, $s\in [N]$,
 is constructed from some graph
$e_s\in \widehat{\HoLB}_{0,1}^{\atop \bigstar\circlearrowright}(m_s,n_s)$
as follows:
\Bi
\item[(i)]
draw new $k$ big white vertices labelled from 1 to k (the ``inputs" of $\hat{e}_s$) and one extra output big white vertex,
\item[(ii)] symmetrize all $m_s$ outputs legs of $e_s$ (if $m_s\geq 1$) and attach them to the unique output white vertex; if $m_s=0$, the output big white vertex receives no incoming edges;
\item[(iii)] partition the set $[n_s]$ if input legs of $e_s$ into $k$ ordered disjoint  subsets
    $$
    [n_s]=I_1\sqcup \ldots \sqcup I_k, \ \ \ \ \#I_i\geq 0, i\in [k],
    $$
and then symmetrize the legs in each subset $I_i$ and attach them (if any) to the $i$-labelled input white vertex.
\Ei

For example, the element
$$
e=\Ba{c}\resizebox{12mm}{!}{
\xy
(0,0)*{\circ}="o",
(-2,5)*{}="2",
(4,5)*{\circ}="3",
(4,10)*{}="u",
(4,0)*{}="d1",
(7,0)*{}="d2",
(10,0)*{}="d3",
(-1.5,-5)*{}="5",
(1.5,-5)*{}="6",
(4,-5)*{}="7",
(-4,-5)*{}="8",
(-2,7)*{_1},
(4,12)*{_2},
(-1.5,-7)*{_2},
(1.5,-7)*{_3},
(10.4,-1.6)*{_6},
(-4,-7)*{_1},
(4,-1.6)*{_4},
(7,-1.6)*{_5},
\ar @{-} "o";"2" <0pt>
\ar @{-} "o";"3" <0pt>
\ar @{-} "o";"5" <0pt>
\ar @{-} "o";"6" <0pt>
\ar @{-} "o";"8" <0pt>
\ar @{-} "3";"u" <0pt>
\ar @{-} "3";"d1" <0pt>
\ar @{-} "3";"d2" <0pt>
\ar @{-} "3";"d3" <0pt>
\endxy}\Ea
 \in \widehat{\HoLB}_{0,1}^{\atop \bigstar\circlearrowright}(2,6)
$$
can generate the following element
$$
\hat{e}=\Ba{c}\resizebox{18mm}{!}{ \xy
(-1.5,5)*{}="1",
(1.5,5)*{}="2",
(9,5)*{}="3",
 (0,0)*{\circ}="A";
  (9,3)*{\circ}="O";
   (5,12)*+{\hspace{2mm}}*\frm{o}="X";
 (-6,-10)*+{_1}*\frm{o}="B";
  (6,-10)*+{_2}*\frm{o}="C";
   (14,-10)*+{_3}*\frm{o}="D";
    (22,-10)*+{_4}*\frm{o}="E";
 "A"; "B" **\crv{(-5,-0)}; 
  "A"; "D" **\crv{(5,-0.5)};
  "A"; "C" **\crv{(-5,-7)};
   "A"; "O" **\crv{(5,5)};
\ar @{-} "O";"C" <0pt>
\ar @{-} "O";"D" <0pt>
\ar @{-} "O";"X" <0pt>
\ar @{-} "A";"X" <0pt>
\ar @{-} "O";"B" <0pt>
 \endxy}
 \Ea \in \f(\widehat{\HoLB}_{0,1}^{\atop \bigstar\circlearrowright})(4)
$$
in the associated polydifferential operad. If we erase the top big white vertex and its all attached edges, then we get from elements of $\f(\widehat{\HoLB}_{0,1}^{\atop \bigstar\circlearrowright})$ precisely
M.\ Kontsevich graphs from \cite{Ko2}. The operad $\f(\widehat{\HoLB}_{0,1}^{\atop \bigstar\circlearrowright})$ admits a filtration by the number of small white vertices (that is, by the number
of vertices coming from the underlying generators of $\widehat{\HoLB}_{0,1}^{\atop \bigstar\circlearrowright}$) which we call from now on {\em internal vertices}.
The big white vertices of graphs from $\f(\widehat{\HoLB}_{0,1}^{\atop \bigstar\circlearrowright})$ are called the {\em external}\, ones. Note that incoming external vertices are {\em not}\, ordered from left to right as one might infer from the pictures above --- they are only labelled by distinct integers.

 Note also that elements of $\widehat{\HoLB}_{0,1}^{\atop \bigstar\circlearrowright}$ may contain elements with no internal vertices at all, for example,
 $$
 \Ba{c}\resizebox{9mm}{!}{ \xy
   (0,9)*+{\hspace{2mm}}*\frm{o}="X";
 (-5,0)*+{_1}*\frm{o}="B";
  (5,0)*+{_2}*\frm{o}="C";
 \endxy
 }\Ea
 \in \widehat{\HoLB}_{0,1}^{\atop \bigstar\circlearrowright}(2).
 $$
 The latter graph admits an automorphism which swaps numerical labels of vertices (cf.\ \cite{MW2,MW3}) and controls the canonical graded commutative multiplication in $\widehat{\odot^\bu}V$. For any $i\in [n]$ the operadic composition
 $$
 \Ba{rccc}
 \circ_i: &\f(\widehat{\HoLB}_{0,1}^{\atop \bigstar\circlearrowright})(n)\ot \f(\widehat{\HoLB}_{0,1}^{\atop \bigstar\circlearrowright})(m) & \lon &
 \f(\widehat{\HoLB}_{0,1}^{\atop \bigstar\circlearrowright})(m+n-1)\\
 & \Ga_1\ot \Ga_2 &\lon & \Ga_1\circ_i \Ga_2
 \Ea
$$
is defined by substituting the graph $\Ga_2$ (with the output external vertex erased so that all edges, if any, connected to that external vertex become ``dangling in the air")) inside the big circle of the $i$-labelled external vertex of $\Ga_1$ and erasing that big circle (so that all edges of $\Ga_1$ connected to the $i$-th external vertex, if any, also become ``dangling in the air"),   and then taking the sum over all possible ways to do the following operations
\Bi
\item[(i)] glueing some (or all or none) hanging edges of $\Ga_2$ to some hanging edges of $\Ga_1$,
\item[(ii)] attaching some (or all or none) hanging edges of $\Ga_2$ to the output external vertex of $\Ga_1$,
\item[(iii)] attaching some (or all or none) hanging edges of $\Ga_1$ to the external input vertices of $\Ga_2$,
\Ei
in such a way that no dangling edges are left. We refer to \cite{MW2,MW3} for concrete examples.

\subsection{Kontsevich formality map as a morphism of dg operads}
M.\ Kontsevich formality map from \cite{Ko2} provides us with a universal quantization of arbitrary (formal) graded Poisson structures. It can understood as a morphism of dg props\footnote{Similarly, a universal formality map behind Drinfeld's deformation quantizations of Lie bialgebras can be understood as a morphism of dg props, see  \cite{MW3}.},
\Beq\label{5: quant map F}
\cF: c\cA ss_\infty \lon \f(\swHoLB_{0,1})
\Eeq
satisfying a certain non-triviality condition (which is given explicitly below). Here $c\cA ss_\infty$ is a  dg operad of {\em curved
$A_\infty$-algebras} defined as the free operad generated by the $\bS$-module

$$
E(n):=\K[\bS_n][n-2]= \mbox{span}
\left(\Ba{c}\resizebox{20mm}{!}{
\begin{xy}
 <0mm,0mm>*{\bullet};<0mm,0mm>*{}**@{},
 <0mm,0mm>*{};<-8mm,-5mm>*{}**@{-},
 <0mm,0mm>*{};<-4.5mm,-5mm>*{}**@{-},
 <0mm,0mm>*{};<0mm,-4mm>*{\ldots}**@{},
 <0mm,0mm>*{};<4.5mm,-5mm>*{}**@{-},
 <0mm,0mm>*{};<8mm,-5mm>*{}**@{-},
   <0mm,0mm>*{};<-11mm,-7.9mm>*{^{\sigma(1)}}**@{},
   <0mm,0mm>*{};<-4mm,-7.9mm>*{^{\sigma(2)}}**@{},
   <0mm,0mm>*{};<10.0mm,-7.9mm>*{^{\sigma(n)}}**@{},
 <0mm,0mm>*{};<0mm,5mm>*{}**@{-},
 \end{xy}}\Ea\ \
\right)_{\sigma\in \bS_n},\ \ \ \forall\ n\geq 0
$$
and equipped with the differential given on the generators by the formula
$$
\delta
\Ba{c}\resizebox{10mm}{!}{
\begin{xy}
<0mm,0mm>*{\bullet},
<0mm,5mm>*{}**@{-},
<-5mm,-5mm>*{}**@{-},
<-2mm,-5mm>*{}**@{-},
<2mm,-5mm>*{}**@{-},
<5mm,-5mm>*{}**@{-},
<0mm,-7mm>*{_{1\ \ \ \ldots\ \ \ n}},
\end{xy}}\Ea
=\sum_{k=0}^{n}\sum_{l=0}^{n-k}
(-1)^{k+l(n-k-l)+1}
\Ba{c}\resizebox{30mm}{!}{
\begin{xy}
<0mm,0mm>*{\bullet},
<0mm,5mm>*{}**@{-},
<4mm,-7mm>*{^{1\ \ \dots \ \ k\qquad\ \    (k+l+1)\ \ \dots\ n}},
<-14mm,-5mm>*{}**@{-},
<-6mm,-5mm>*{}**@{-},
<20mm,-5mm>*{}**@{-},
<8mm,-5mm>*{}**@{-},
<0mm,-5mm>*{}**@{-},
<0mm,-5mm>*{\bullet};
<-5mm,-10mm>*{}**@{-},
<-2mm,-10mm>*{}**@{-},
<2mm,-10mm>*{}**@{-},
<5mm,-10mm>*{}**@{-},
<0mm,-12mm>*{_{k+1\ \dots\ k+l }},
\end{xy}}\Ea.
$$
It is non-cofibrant and acyclic (as the dg operad $\f(\widehat{\HoLB}_{0,1}^{\atop \bigstar\circlearrowright})$).

\sip

The non-triviality condition on the map (\ref{5: quant map F}) reads as the following approximation on the values of $\cF$ on the generating $n$ corollas
of $c\cA ss_\infty$ for any $n\geq 0$
 (modulo graphs with the number of internal vertices $\geq 2$ whose linear span is denoted below by $O(2)$)
\Beq\label{5: Boundary cond for formality map}
\cF\left(
\Ba{c}\resizebox{10mm}{!}{
\begin{xy}
<0mm,0mm>*{\bullet},
<0mm,5mm>*{}**@{-},
<-5mm,-5mm>*{}**@{-},
<-2mm,-5mm>*{}**@{-},
<2mm,-5mm>*{}**@{-},
<5mm,-5mm>*{}**@{-},
<0mm,-7mm>*{_{1\ \ \ \ldots\ \ \ n}},
\end{xy}}\Ea
\right)=\left\{
\Ba{ll}
 \Ba{c}\resizebox{9mm}{!}{ \xy
   (0,9)*+{\hspace{2mm}}*\frm{o}="X";
 (-5,0)*+{_1}*\frm{o}="B";
  (5,0)*+{_2}*\frm{o}="C";
 \endxy}\Ea
 + \sum_{p\geq 0}\frac{1}{p!}
 \Ba{c}\resizebox{10mm}{!}{ \xy
 (0,2.8)*{^p};
 (0,1)*{...};
 (0,-3)*{\circ}="a";
   (0,9)*+{\hspace{2mm}}*\frm{o}="X";
 (-7,-9)*+{_1}*\frm{o}="B";
  (7,-9)*+{_2}*\frm{o}="C";
   "a"; "X" **\crv{(-5,-0)}; 
   "a"; "X" **\crv{(+5,-0)}; 
   "a"; "X" **\crv{(9,-0)}; 
   "a"; "X" **\crv{(-9,-0)}; 
  \ar @{-} "a";"B" <0pt>
  \ar @{-} "a";"C" <0pt>
 \endxy}\Ea + O(2) & \text{if}\ n=2\\
\sum_{p\geq 0}\frac{1}{p!}
 \Ba{c}\resizebox{13mm}{!}{ \xy
 (0,2.8)*{^p};
 (0,1)*{...};
 (3.5,-10)*{...};
 (0,-3)*{\circ}="a";
   (0,9)*+{\hspace{2mm}}*\frm{o}="X";
 (-10,-10)*+{_1}*\frm{o}="B";
  (-3,-10)*+{_2}*\frm{o}="C";
  (10,-10)*+{_{k}}*\frm{o}="E";
   "a"; "X" **\crv{(-5,-0)}; 
   "a"; "X" **\crv{(+5,-0)}; 
   "a"; "X" **\crv{(9,-0)}; 
   "a"; "X" **\crv{(-9,-0)}; 
  \ar @{-} "a";"B" <0pt>
  \ar @{-} "a";"C" <0pt>
  \ar @{-} "a";"E" <0pt>
 \endxy}\Ea + O(2) & \text{otherwise}\\
 \Ea
\right.
\Eeq
where the summations $\sum_{p\geq 0}$
 run over the number of edges connecting the internal vertex to the external out-vertex. A morphism of dg operads (\ref{5: quant map F}) satisfying the above non-triviality condition is called a {\em formality map} (after \cite{Ko2}).

\subsection{Deformation complexes of morphisms of props}
Let $\cP$ be an arbitrary dg free prop, and $\cQ$ an arbitrary dg
prop, and $f: \cP \rar \cQ$ a morphism between them. Then there is a standard construction of the {\em deformation complex}\,
$\Def(\cP \stackrel{f}{\rar}\cQ)$ of the morphism $f$ described in several ways in \cite{MV}; in general, $\Def(\cP \stackrel{f}{\rar}\cQ)$ is a filtered $\caL ie_\infty$ algebra.
 This construction builds on earlier works
which describe deformation complexes of morphisms of dg {\em operads} \cite{KS,VdL}.
The constructions in \cite{MV} generalize straightforwardly to the case when
$\cP$ and $\cQ$ are dg {\em wheeled}\, props.
For example,
when $\cP=\cQ=\wh{\HoLB}_{c,d}^{\atop \bigstar\circlearrowright}$ and $f$ is the identity map, then the associated deformation complex
$$
\Def(\wh{\HoLB}_{c,d}^{\atop \bigstar\circlearrowright}\stackrel{\Id}{\lon}
\wh{\HoLB}_{c,d}^{\atop \bigstar\circlearrowright})[1]\simeq \Der({\HoLB}_{c,d}^{\atop \bigstar\circlearrowright})
$$
is, up to the degree shift, is precisely the derivation complex of
$\wh{\HoLB}_{c,d}^{\atop \bigstar\circlearrowright}$ (but the Lie algebra structure is different!). The machinery of \cite{KS,MV,VdL} gives us  a well-defined
dg Lie algebra
$$
\Def\left(c\cA ss_\infty \stackrel{\cF}{\lon} \f(\wh{\HoLB}_{0,1}^{\atop \bigstar\circlearrowright})\right)
$$
 which controls the deformation theory of any formality map $\cF$. Our second main result in this paper is the computation of its cohomology in terms of the M.\ Kontsevich graph complex $\mathsf{fGC}^{\geq 2}_2$.

\subsection{Theorem (Classification of formality maps)}\label{5: Corollary on GCor and Def(assb to Dlie)} {\em For any formality morphism $\cF$
$$
\cF:
c\Ass_\infty {\lon} \f(\wh{\HoLB}_{0,1}^{\atop \bigstar\circlearrowright})
$$
there is a canonically associated morphism
 of complexes
 $$
f_\cF:  \mathsf{fGC}^{\geq 2}_2 \lon \Def\left(c\Ass_\infty \stackrel{\cF}{\lon} \f(\wh{\HoLB}_{0,1}^{\atop \bigstar\circlearrowright})\right)[1]
 $$
 which is a quasi-isomorphism.}
\begin{proof} The proof of this Theorem is very similar to the proof of Proposition 5.4.1 in \cite{MW3} and is based essentially on the contractibility of the permutahedra polytopes. Let us first explain the naturality of the morphism $f_{\cF}$.
Any derivation of the dg wheeled prop $\wh{\HoLB}_{1,0}^{\atop \bigstar\circlearrowright}$, that is, any deformation $D$ of the identity automorphism
of $\wh{\HoLB}_{0,1}^{\atop \bigstar\circlearrowright}$,
$$
D\in \Def\left( \wh{\HoLB}_{0,1}^{\atop \bigstar\circlearrowright}\stackrel{\Id}{\lon} \wh{\HoLB}_{1,0}^{\atop \bigstar\circlearrowright}\right)
$$
induces an associated deformation of the identity automorphism of $\f(\wh{\HoLB}_{0,1}^{\atop \bigstar\circlearrowright})$,
$$
D\in \Def\left(\f(\wh{\HoLB}_{1,0}^{\atop \bigstar\circlearrowright})\stackrel{\Id}{\lon}
\f(\wh{\HoLB}_{1,0}^{\atop \bigstar\circlearrowright})
\right)
$$
and hence, via the composition of $D$ with the given map $\cF$, gives us a canonical morphism of complexes
$$
g_{\cF}: \Def\left(\wh{\HoLB}_{0,1}^{\atop \bigstar\circlearrowright}\stackrel{\Id}{\lon}  \wh{\HoLB}_{0,1}^{\atop \bigstar\circlearrowright}\right)\lon
\Def\left(c\Ass_\infty \stackrel{\cF}{\lon} \f( \wh{\HoLB}_{0,1}^{\atop \bigstar\circlearrowright})\right)
$$
or, equivalently,
\Beq\label{5: g_cF}
g_{\cF}: \Der(\wh{\HoLB}_{0,1}^{\atop \bigstar\circlearrowright}) \lon
\Def\left(c\Ass_\infty \stackrel{\cF}{\lon} \f( \wh{\HoLB}_{0,1}^{\atop \bigstar\circlearrowright})\right)[1]
\Eeq
Composing this map $g_{\cF}$ with the canonical quasi-isomorphism $F$ from Proposition 1.1.1, we obtain the required map $f_{\cF}$. Thus to prove the theorem it is enough to prove that the map $g_{\cF}$ is a quasi-isomorphism. Which is easy.

\sip

 Both complexes in (\ref{5: g_cF}) admits filtrations by the number of edges in the graphs, and the map $g_{\cF}$ preserves these filtrations, and hence
induces a morphism of the associated spectral sequences,
$$
g^r_{\cF}: (\cE_r\Der(\wh{\HoLB}_{0,1}^{\atop \bigstar\circlearrowright}), d_r) \lon \left(\cE_r
\Def\left(c\Ass_\infty \stackrel{\cF}{\lon} \f( \wh{\HoLB}_{0,1}^{\atop \bigstar\circlearrowright})\right)[1], \delta_r\right).
$$
The induced differential $d_0$ on the initial page of the spectral sequence of the l.h.s.\ is trivial, $d_0=0$. The induced differential  on the initial page of the spectral sequence of the r.h.s. is not trivial and is determined by the
following summand in $\cF$ (see (\ref{5: Boundary cond for formality map})),
$$
 \Ba{c}\resizebox{9mm}{!}{ \xy
   (0,9)*+{\hspace{2mm}}*\frm{o}="X";
 (-5,0)*+{_1}*\frm{o}="B";
  (5,0)*+{_2}*\frm{o}="C";
 \endxy}\Ea
$$
Hence the differential $\delta_0$ acts only on big input  white vertices of graphs from $\Def\left(c\Ass_\infty \stackrel{\cF}{\rar} \f( \wh{\HoLB}_{0,1}^{\atop \bigstar\circlearrowright})\right)[1] $ by splitting each such big white vertex
$
 \Ba{c}\resizebox{4mm}{!}{
 \xy
 (0,1.5)*+{_v}*\frm{o};
 \endxy
 }\Ea
$
into two big white vertices
$
 \Ba{c}\resizebox{10mm}{!}{ \xy
 (-5,2)*+{_{v'}}*\frm{o}="B";
  (5,2)*+{_{v''}}*\frm{o}="C";
 \endxy}\Ea
$
and redistributing all edges (if any) attached to $v$  in all possible ways among the new vertices $v'$ and $v''$. The cohomology
$$ \cE_1
\Def\left(c\Ass_\infty \stackrel{\cF}{\lon} \f( \wh{\HoLB}_{0,1}^{\atop \bigstar\circlearrowright})\right)[1]= H\left(\cE_0
\Def\left(c\Ass_\infty \stackrel{\cF}{\lon} \f( \wh{\HoLB}_{0,1}^{\atop \bigstar\circlearrowright})\right)[1], \delta_0\right)$$
is spanned by graphs all of whose white vertices are precisely univalent and skew symmetrized (see, e.g., Theorem 3.2.4 in \cite{Me-p} where this result is obtained from the cell complexes of permutahedra, or Appendix A in \cite{Wi1} for another purely algebraic argument) and hence is isomorphic (after erasing these no more needed big white vertices)  to $\Der(\swHoLB_{0,1})$ as a graded vector space.
The boundary  condition (\ref{5: Boundary cond for formality map}) says that the induced differential $\delta_1$ in the complex
$
 \cE_1
\Def\left(c\Ass_\infty \stackrel{\cF}{\lon} \f( \wh{\HoLB}_{0,1}^{\atop \bigstar\circlearrowright})\right)[1]$
agrees precisely with the induced differential $d_1$ in
$\cE_1\Der(\swHoLB_{0,1})$ so that the induced morphism of the next pages of the spectral sequences,
$$
g_{\cF}^1: ( \wh{\HoLB}_{0,1}^{\atop \bigstar\circlearrowright}), d_1)\lon
\left(\cE_1\Def(c\cA ss_\infty \stackrel{\cF}{\rar} \f(\wh{\HoLB}_{0,1}^{\atop \bigstar\circlearrowright}), \delta_1\right)
$$
is an isomorphism. By the spectral sequence comparison theorem, the morphism $g_{\cF}$ is a quasi-isomorphism.
\end{proof}

\mip
We conclude that for any $i\in \Z$,
$$
H^{i+1}\left(\Def\left(c\cA ss_\infty \stackrel{\cF}{\rar} \f(\wh{\HoLB}_{0,1}^{\atop \bigstar\circlearrowright})\right)\right)=H^i(\fGC_2^{\geq 2})
$$
The special case $i=0$ reads as
$$
H^1\left(\Def\left(c\cA ss_\infty \stackrel{\cF}{\rar} \f(\wh{\HoLB}_{0,1}^{\atop \bigstar\circlearrowright})\right)\right)=H^0(\fGC_2^{\geq 2}) =\grt
$$
which is equivalent to the main result of the remarkable paper \cite{Do} by V.\ Dolgushev.

\sip

\bip

\def\cprime{$'$}

\end{document}